\documentclass[11pt]{article}
\usepackage{amsmath} 
\usepackage{amsthm} 
\usepackage{amssymb} 
\usepackage{hyperref}
\usepackage{comment} 
\usepackage{enumitem} 
\usepackage{color}

\newtheorem{thm}{Theorem}
\newtheorem{inspr}[thm]{}

\newenvironment{definitie}{\begin{itemize}\item[ ]\hspace{-26pt}\bf Definition \rm }{\end{itemize}}
\newenvironment{notatie}{\begin{itemize}\item[ ]\hspace{-26pt}\bf Notation \rm }{\end{itemize}}
\newenvironment{voorbeeld}{\begin{itemize}\item[ ]\hspace{-26pt}\bf Example \rm }{\end{itemize}}
\newenvironment{stelling}{\begin{itemize}\item[ ]\hspace{-26pt}\bf Theorem \rm }{\end{itemize}}
\newenvironment{propositie}{\begin{itemize}\item[ ]\hspace{-26pt}\bf Proposition \rm }{\end{itemize}}
\newenvironment{lemma}{\begin{itemize}\item[ ]\hspace{-26pt}\bf Lemma \rm }{\end{itemize}}
\newenvironment{opmerking}{\begin{itemize}\item[ ]\hspace{-26pt}\bf Remark \rm }{\end{itemize}}
\newenvironment{voorwaarde}{\begin{itemize}\item[ ]\hspace{-26pt}\bf Condition \rm }{\end{itemize}}
 
\newcommand{\defin}{\begin{inspr}\begin{definitie}}  
\newcommand{\edefin}{\end{definitie}\end{inspr}}
\newcommand{\notat}{\begin{inspr}\begin{notatie}}  
\newcommand{\enotat}{\end{notatie}\end{inspr}}
\newcommand{\voorb}{\begin{inspr}\begin{voorbeeld}}
\newcommand{\evoorb}{\end{voorbeeld}\end{inspr}}
\newcommand{\stel}{\begin{inspr}\begin{stelling}}
\newcommand{\estel}{\end{stelling}\end{inspr}}
\newcommand{\prop}{\begin{inspr}\begin{propositie}}
\newcommand{\eprop}{\end{propositie}\end{inspr}}
\newcommand{\lem}{\begin{inspr}\begin{lemma}}
\newcommand{\elem}{\end{lemma}\end{inspr}}
\newcommand{\opm}{\begin{inspr}\begin{opmerking}}
\newcommand{\eopm}{\end{opmerking}\end{inspr}}
\newcommand{\voorw}{\begin{inspr}\begin{voorwaarde}}
\newcommand{\evoorw}{\end{voorwaarde}\end{inspr}}
\newcommand{\bew}{\vspace{-0.3cm}\begin{itemize}\item[ ] \bf Proof\rm: }
\newcommand{\ebew}{\hfill $\qed$ \end{itemize}}
\newcommand{\ssnl}{\vskip 3pt} 
\newcommand{\snl}{\vskip 7pt} 
\newcommand{\nl}{\vskip 12pt} 

\newcommand{\ot}{\otimes}

\newcommand{\inv}{^{-1}}

\newcommand{\tl}{\triangleleft}
\newcommand{\tr}{\triangleright}
\newcommand{\tussenen}{\qquad\quad\text{and}\qquad\quad}

\newcommand{\blauw}{\color{blue}}

\newcommand{\Chi}{\raisebox{2pt}{\ensuremath{\chi}}}

\numberwithin{thm}{section}   
\numberwithin{equation}{section} 
\begin{document}

\centerline{\bf \Large Finite quantum hypergroups}
\vspace{13pt}
\centerline{\it Magnus B. Landstad \rm $^{(1)}$ and \it Alfons Van Daele \rm $^{(2)}$}
\bigskip\bigskip

{\bf Abstract}
\nl
A \emph{finite quantum hypergroup}
is a finite-dimensional unital algebra $A$ over the field of complex numbers. There is a coproduct on $A$, a coassociative map from $A$ to $A\ot A$ assumed to be unital, but it is not required to be an algebra homomorphism. There is a counit  that is supposed to be a homomorphism. Finally, the main extra requirement is the existence of a faithful left integral with the right properties.
\ssnl
For such a finite quantum hypergroup, the dual can be constructed. It is again a finite quantum hypergroup.
\ssnl
The more general concept of an \emph{algebraic quantum hypergroup} is studied in \cite{De-VD1, De-VD2}. If the underlying algebra of an algebraic quantum hypergroup is finite-dimensional, it is a finite quantum hypergroup in the sense of this paper. 
\ssnl
Here we treat  finite quantum hypergroups independently with an emphasis on the development of the notion. It is meant to clarify the various steps taken in \cite{La-VD3a} and \cite{La-VD3b}. In \cite{La-VD3b} we introduce and study a still more general concept of quantum hypergroups. It not only contains the algebraic quantum hypergroups from \cite{De-VD1}, but also some topological cases. We naturally encounter these quantum hypergroups in our work on bicrossproducts, see \cite{La-VD4, La-VD5, La-VD7}.
\ssnl
We include some  examples to illustrate the theory. One kind is coming from a finite group and a subgroup. The other examples are taken from the bicrossproduct theory.
\nl
Date: {\it 26 September 2022} 
\nl
\vskip 1cm
\hrule
\medskip
\begin{itemize}
\item[$^{(1)}$] Department of Mathematical Sciences, Norwegian University of Science  and  Technology, NO-7491 Trondheim, Norway. E-mail: \texttt{magnus.landstad@ntnu.no}
\item[$^{(2)}$] Department of Mathematics, University of Leuven, Celestijnenlaan 200B,
B-3001 Heverlee, Belgium. E-mail: \texttt{alfons.vandaele@kuleuven.be}
\end{itemize}
\newpage

\setcounter{section}{-1}  

\section{\hspace{-17pt}. Introduction} \label{s:introduction}

Let $A$ and $B$ be associative algebras over the field $\mathbb C$ of complex numbers. In the first place, we do not require these algebras to be unital but we do need that the product is non-degenerate. This means that, given $a\in A$, then $a=0$ if $aa'=0$ for all $a'\in A$ and similarly on the other side. 
Further we assume that we have a  pairing of $A$ with $B$. This is a bilinear map $(a,b)\mapsto \langle a,b\rangle$ from $A\times B$ to $\mathbb C$. It is assumed to be non-degenerate. So, given $a\in A$, then $a=0$ if $\langle a,b\rangle=0$ for all $b\in B$ and similarly, given $b\in B$, we must have $b=0$ if $\langle a,b\rangle=0$ for all $a\in A$.
\ssnl
This will be the core of the treatment of quantum hypergroups in \cite{La-VD3b}. In this note, we will only work with \emph{finite-dimensional unital algebras}. Then the product is automatically non-degenerate.
In \cite{La-VD3b} we will treat the more general case of possibly infinite-dimensional, non-unital algebras. See also \cite{La-VD3a} for a still more general duality concept.
\nl 
\bf Pairing of finite-dimensional unital algebras  \rm
\snl
Now let $A$ and $B$ be \emph{finite-dimensional unital} algebras over $\mathbb C$. We use $1$ for the unit in $A$ and in $B$. Assume that we have a pairing of $A$ with $B$ as above.
\ssnl
The product in $B$ induces a coproduct $\Delta_A$ on $A$ and the product in $A$ gives a coproduct $\Delta_B$ on $B$ by the formulas 
\begin{equation*}
\langle a,bb'\rangle=\langle \Delta_A(a),b\ot b'\rangle
\tussenen
\langle aa',b\rangle=\langle a\ot a', \Delta_B(b)\rangle
\end{equation*}
where $a,a'\in A$ and $b,b'\in B$. In what follows, we will drop the indices and use $\Delta$ both for $\Delta_A$ and $\Delta_B$. The coproduct on $A$ is then a coassociative linear map $\Delta: A\mapsto A\ot A$ in the sense that 
$(\Delta\ot\iota)\Delta=(\iota\ot\Delta)\Delta$ where $\iota$ is the identity map on $A$. Similarly for the coproduct on $B$.
\ssnl
The identity in $B$ induces a linear functional $\varepsilon_A$ on $A$ and the identity in $A$ gives a linear functional $\varepsilon_B$ on $B$ by the formulas
\begin{equation*}
\varepsilon_A(a)=\langle a,1\rangle
\tussenen
\varepsilon_B(b)=\langle 1,b\rangle
\end{equation*}
where $a\in A$ and $b\in B$. Again we drop the indices and we will use $\varepsilon$ for $\varepsilon_A$ and $\varepsilon_B$. The map $\varepsilon$ on $A$ is a counit in the sense that 
\begin{equation*}
(\varepsilon\ot\iota)\Delta(a)=a
\tussenen
(\iota\ot\varepsilon)\Delta(a)=a
\end{equation*}
for all $a\in A$. Similarly for $\varepsilon$ on $B$.
\ssnl
Dual pairs of algebras like above appear in various cases, with different extra properties of the coproducts and the counits.
\ssnl
If  $A$ is a finite-dimensional \emph{Hopf algebra} and $B$ its dual Hopf algebra, then the coproduct on $A$ is a unital homomorphism  and the counit on $A$ is a homomorphism. Similarly for the coproduct and the counit on $B$. In particular $\Delta(1)=1\ot 1$, both for $A$ and $B$. This is equivalent with the counits being homomorphisms.
\ssnl
If $A$ is a \emph{weak Hopf algebra}, we still have that $\Delta$ is a homomorphism, but it is no longer assumed to be unital. The counit is no longer  a homomorphism. 
\ssnl
On the other hand, for a pairing of \emph{finite quantum hypergroups}, as in \cite{De-VD1}, we no longer have that $\Delta$ is a homomorphism, but keep the condition  $\Delta(1)=1\ot 1$, both on $A$ and $B$. 
\ssnl
Loosely speaking, for a quantum group, the coproduct is a unital homomorphism. For a quantum groupoid it is no longer required to be unital, while for a quantum hypergroup it is no longer a homomorphism. 
\ssnl
For precise definitions, we refer to the item `Basic References' further in this introduction. The notion of a finite quantum group as treated in this paper is found in Definition \ref{defin:3.10} in Section \ref{s:finite}.
\nl
\bf Pairing of $^*$-algebras \rm
\snl
Now let $A$ and $B$ be unital $^*$-algebras, still finite-dimensional, with a non-degenerate pairing.  
\ssnl
Then there are linear maps $S_A$ on $A$ and $S_B$ on $B$ determined by
\begin{equation*}
\langle a,b^*\rangle=\overline{\langle S_A(a)^*,b\rangle}
\tussenen
\langle a^*,b\rangle=\overline{\langle a,S_B(b)^*\rangle} 
\end{equation*}
for all $a\in A$ and $b\in B$. The maps $S_A$ and $S_B$ are each others adjoints in the sense that $\langle S_A(a),b\rangle=\langle a,S_B(b)\rangle$ for all $a,b$. See Proposition \ref{prop:1.10} in Section \ref{s:pairs}. Also here  we further drop the indices and write $S$ in both cases.
\ssnl
In the case of $^*$-algebras, it is natural to assume that the induced coproduct $\Delta$ on $A$ is a $^*$-map, i.e.\ that $\Delta(a^*)=\Delta(a)^*$ for all $a$. This is equivalent with saying that the map $b\mapsto S(b)^*$ is an algebra map so that $S$ is an anti-isomorphism of $B$. If this is also assumed for $\Delta$ on $B$, it follows that the map $S$ is an anti-isomorphism of $A$ that flips the coproduct. Similarly on $B$.
\ssnl
In the case of a Hopf $^*$-algebra, we take for $S$ the antipode on $A$ and make the dual Hopf algebra $B$  into a $^*$-algebra using the formula $\langle a,b^*\rangle=\overline{\langle S(a)^*,b\rangle}$. Then $B$ is again a Hopf $^*$-algebra. It is used that $S$ is an anti-isomorphism that flips the coproduct.
\ssnl
For a dual pair of weak Hopf $^*$-algebras, this is still true, as well as for a pair of $^*$-algebraic quantum hypergroups.
\ssnl
We will not restrict ourselves to $^*$-algebras because the involution is not essential for our theory. On the other hand, $^*$-algebras are used to motivate the treatment of the antipode in the general theory. This is possible because, as we showed above, the antipodes are determined by the pairing.
\nl
\bf The antipode for pairs of algebras, not necessarily $^*$-algebras \rm
\snl
For a pairing of Hopf algebras, we also have  antipodes. As antipodes are unique, they can be recovered from the pairing of the algebras. They are anti-isomorphisms of the algebras, flipping the coproducts. However, this property does not completely determine the antipodes. If $A$ and $B$ are abelian, also the identity map will satisfy this property.
\ssnl
This means that we must use another way to obtain the antipodes from the pairing of algebras. It can be done in various ways. Using the Sweedler notation we can write
\begin{equation*}
a\ot 1=\sum_{(a)}\Delta(a_{(1)})(1\ot S(a_{(2)}))=\sum_i \Delta(p_i)(1\ot q_i)
\end{equation*}
and then we have 
\begin{equation*}
\sum_i (1\ot p_i)\Delta(q_i)=\sum_{(a)}(1\ot a_{(1)})\Delta(S(a_{(2)}))=S(a)\ot 1.
\end{equation*}
This can be used to recover the antipode from the coproduct.
\ssnl
We have the same property for a multiplier Hopf algebra with integrals, see \cite{VD-mha, VD-alg}. Also here all information can be recovered from the algebra $A$, its dual algebra $\widehat A$ and the pairing between the two. In fact there are plenty of other examples with the same phenomenon. 
\ssnl
This also applies to the pairing of algebraic quantum hypergroups as studied in \cite{De-VD1}. Here, in this paper, we only consider  the finite-dimensional case. We give a treatment  that also works  for the quantum hypergroups we will be studying in \cite{La-VD3b}. 
\ssnl
More precisely, in this paper, we find conditions on a pair of finite-dimensional unital algebras to be a pair of finite quantum hypergroups. 
\ssnl
In \cite{T-VD-W} the duality of algebraic quantum groupoids is treated from the same point of view. But that case is very different from, and more involved than the one studied here. 
\nl
\bf Style of the paper \rm
\nl
The main intention of writing this note is to provide the underlying ideas for the development of a more general theory of quantum hypergroups in \cite{La-VD3a} and \cite{La-VD3b}. The finite-dimensional case is  useful for 
understanding the purely algebraic aspects of the infinite-dimensional theory. This is also the case here.
\ssnl
Finite quantum hypergroups are  special cases of algebraic quantum hypergroups (as studied in \cite{De-VD1} and \cite{De-VD2}). On the other hand, the quantum hypergroups in \cite{La-VD3b} are more general than in \cite{De-VD1}. The paper \cite{De-VD1} could also be used for motivating the approach to quantum hypergroups in \cite{La-VD3b} but still, we think the finite-dimensional case is more suited for this purpose.
\ssnl
Furthermore, the treatment of finite quantum hypergroups as in this paper is independent from the earlier works and the approach is slightly different. We give more details and arguments that are behind this approach. We are convinced that this will contribute to a better understanding of the theory of quantum hypergroups in general.
\ssnl
Finally, we have included examples to illustrate our notions. They are mostly new and there is little overlap with examples of algebraic quantum hypergroups already given in the earlier papers.

\nl
\bf Content of the paper \rm 
\nl
In \emph{Section} \ref{s:pairs} we consider a pair of \emph{finite-dimensional unital algebras} $A$ and $B$ with a non-degenerate pairing. As we discussed already, the product in one algebra determines a coproduct on the other one and the identity in one algebra gives a counit on the other one. Various extra properties of these coproducts and counits are considered. We are particularly interested in  the pairing of $^*$-algebras. The main, but natural condition is now that the coproducts are $^*$-maps. Then there are natural candidates for the antipodes. This serves as a motivation for the introduction of antipodes in this context in the next section.  
\ssnl
In \emph{Section} \ref{s:integrals} we introduce the notion of left and right integrals. The notion involves the existence of an antipode. The definition is first motivated by a result about $^*$-algebras because of the special role of the antipode in this case.  For a faithful left integral, the associated antipode is unique. Further, given a faithful left integral on $A$, there exists a faithful right integral on $B$. This extends the duality from the algebras $A$ and $B$ to the case where these algebras have integrals. In \cite{La-VD3a}, we do the same but for possibly infinite-dimensional algebras.
\ssnl
In \emph{Section} \ref{s:finite} we work towards the main definition of this note. We first study the relation between integrals and invariant functionals. An important result is the uniqueness of integrals when the coproduct is unital. This eventually leads to the definition of a finite quantum hypergroup. The duality obtained in the previous section applies here and gives a duality for finite quantum hypergroups.
\ssnl
In the  subsequent sections we illustrate our theory with examples. There are two sets of examples of a different nature. 
\ssnl
In \emph{Section} \ref{s:hecke} we treat the case of a finite group $G$ with a subgroup $H$ and the $^*$-algebra of functions on $G$ with the property that $f(hpk)=f(p)$ for all $p\in G$ and $h,k\in H$. We have a coproduct $\Delta$ on $A$ so that $(A,\Delta)$ is a finite $^*$-quantum hypergroup. We have an explicit construction of the dual. This is the classical Hecke algebra and it is the finite-dimensional motivating example for the study of algebraic quantum hypergroups as in \cite{De-VD1,De-VD2}. 
\ssnl
We make it more concrete with the subgroup generated by a single permutation in the group $S_3$ of permutations of a set with $3$ elements in \emph{Section} \ref{s:two-dim}. We also derive some generalizations here. And we finish the section with related examples of pairs where invariant functionals exist that are not integrals.
\ssnl
In \emph{Section} \ref{s:two-sub} we look at a  group $G$, not necessarily finite,  with two finite subgroups $H$ and $K$ satisfying $H\cap K=\{e\}$.
We associate a dual pair of finite quantum hypergroups to such a pair of subgroups. We show that we have a pair of finite quantum groups (Hopf algebras) if and only if $HK=KH$. This is the finite-dimensional version of the quantum hypergroups we encounter in the theory of bicrossproducts in \cite{La-VD4, La-VD5, La-VD7}. We consider the special case of groups $H$ and $K$ sitting in their free product. Moreover, we illustrate this for groups with only two elements.
\ssnl
We finish with some conclusions and suggestions for further research in Section \ref{s:conclusions}.
\nl
\bf Notations and conventions \rm
\nl
We work with associative algebras  over the field $\mathbb C$ of complex numbers. They are always assumed to be finite-dimensional and unital. 
 \ssnl
We use $1$ for the identity in any algebra we consider. For the identity in a group, we use $e$. Finally $\iota$ is used for the identity map. 
\ssnl
A linear functional $\omega$ on $A$ is called faithful if, given $a'\in A$, we have $a=0$ if either $\omega(aa')=0$ for all $a\in A$ or if $\omega(a'a)=0$ for all $a'\in A$. For a positive linear functional $\omega$ on a $^*$-algebra, this is equivalent with the requirement that $\omega(a^*a)=0$ only if $a=0$. A finite-dimensional $^*$-algebra with a faithful positive linear functional is an operator algebra. 
\ssnl
In this note, a coproduct on an algebra $A$ is  just a coassociative linear map from $A$ to the tensor product $A\ot A$. We call it coabelian if $\zeta\circ\Delta=\Delta$ where $\zeta$ is the linear map on $A\ot A$ defined by $\zeta(a\ot a')=a'\ot a$ for all $a,a'\in A$.
\ssnl
Often in this paper, we will use the same symbol for objects associated with different algebras. We not only use e.g.\ $1$ to denote the identity in all our algebras, we will do the same for coproducts, counits and antipodes.
\ssnl
We are inspired by the theory of locally compact groups and operator algebras. This has some consequences for the terminology, notations and techniques we use. This also explains why we only work with algebras over the complex numbers and why we also treat $^*$-algebras.
\nl
\bf Basic references \rm
\nl
For the theory of Hopf algebras, we refer to the basic works by Abe \cite{Ab} and Sweedler \cite{Sw}. See also \cite{R-bk} for a more recent treatment.
 For multiplier Hopf algebras and algebraic quantum groups, the main references are \cite{VD-mha} and \cite{VD-alg}. 
The original theory of quantum hypergroups is found in \cite{De-VD1} and \cite{De-VD2}. There is also the earlier work on compact quantum hypergroups in \cite{Ch-Va}. The theory of weak Hopf algebras and weak multiplier Hopf algebras is developed in a series of papers. References can e.g.\ be found in \cite{VD-W2}. In particular, the pairing of weak multiplier Hopf algebras is treated in \cite{T-VD-W}. 
\ssnl
In \cite{La-VD3} a shortened version of this note is found.
\nl
\bf Acknowledgments \rm
\nl
We want to thank the Mathematics Departments of NTNU (Norway) and KU Leuven (Belgium) for the opportunities we get to continue our research after our official retirements. 
\nl

\section{\hspace{-17pt}. Dual pairs of algebras and $^*$-algebras} \label{s:pairs} 

We start this section with some simple and well-known properties of finite-dimensional algebras over the field  $\mathbb C$ of complex numbers. Further, we consider a coproduct with a counit. Then the dual space is again an algebra and we get a pairing of algebras. Special attention goes to the case of $^*$-algebras. As explained already, this case is used to motivate the main formulas in the non-involutive case.
\nl
\bf Finite-dimensional algebras with a faithful linear functional \rm
\nl
We start with a remark about faithfulness of a linear functional on an algebra. We have one-sided notions. Left faithfulness of a linear functional $\omega$ on an algebra $A$ means that, given $a\in A$, we have $a=0$ if $\omega(ac)=0$ for all $c\in A$.  Right faithfulness  means that, given $c$, we have $c=0$ if $\omega(ac)=0$ for all $a$. For a finite-dimensional algebra these two properties are equivalent as we show in the following proposition.

\prop
Let $A$ be a finite-dimensional algebra. A linear functional $\omega$ on $A$ is left faithful if and only if it is right faithful. 
\eprop
\bew
i) Assume that $\omega$ is left faithful. Then the map $a\mapsto \omega(a\,\cdot\,)$ is injective. Because $A$ is finite-dimensional, it is also surjective. Now suppose that $c\in A$ and that $\omega(ac)=0$ for all $a$. By the previous result, this means that $\omega'(c)=0$ for all $\omega'$. This implies that $c=0$. Therefore $\omega$ is right faithful. 
\ssnl
ii) The converse property is proven in the same way.
\ebew

\ssnl
Remark that this result is no longer true for infinite-dimensional algebras. In that case, we require both left and right faithfulness to call a functional faithful. 
\nl
We see from the argument in the proof that the following is true.

\prop\label{prop:1.2}
Let $\omega$ be a faithful functional on an  algebra $A$.  For any linear functional $f$ on $A$, there exist elements $c,c'$ in $A$ so that $f(a)=\omega(ac)$ and $f(a)=\omega(c'a)$ for all $a$. 
\eprop

\bew
The maps $c\mapsto \omega(\,\cdot\,c)$ and $c'\mapsto \omega(c'\,\cdot\,)$ are injective because $\omega$ is faithful. Then they are surjective because $A$ is finite-dimensional.
\ebew 
This property has the following consequence. We thank J. Vercruyse for providing us the following simple argument to prove this result.

\prop\label{prop:1.3}
If $A$ is a finite-dimensional algebra with a faithful linear functional, then it is unital.
\eprop
\bew
Let $\omega$ be a faithful linear functional on $A$. By the previous result we have elements $e,f\in A$ satisfying
\begin{equation*}
\omega=\omega(e\,\cdot\,)=\omega(\,\cdot\,f).
\end{equation*}
Then $\omega(ac)=\omega(eac)$ for all $c$ and by the faithfulness of $\omega$ we have $a=ea$. This holds  for all $a$. Similarly we have $a=af$ for all $a$. Then $e=f$ and it is a unit of $A$.
\ebew

\prop\label{prop:1.4}
For any faithful linear functional $\omega$ there exist an automorphism $\sigma$ satisfying $\omega(ac)=\omega(c\sigma(a))$ for all $a,c$. Moreover $\omega(\sigma(a))=\omega(a)$ for all $a$.
\eprop

\bew
Given $a\in A$, the functional $c\mapsto \omega(ac)$ is of the form $c\mapsto \omega(ca')$ for some $a'$ by Proposition \ref{prop:1.2}. Because $\omega$ is faithful, the element $a'$ is uniquely determined by $a$. Denote it by $\sigma(a)$. Then we have
\begin{equation*}
\omega(c\sigma(aa'))=\omega(aa'c)=\omega(a'c\sigma(a))=\omega(c\sigma(a)\sigma(a'))
\end{equation*}
for all $c$ and hence, by the faithfulness of $\omega$, we must have $\sigma(aa')=\sigma(a)\sigma(a')$. The map $\sigma$ is injective and so it is an automorphism of $A$. By Proposition \ref{prop:1.3} the algebra has a unit. Then we can take $c=1$ in the defining formula and we find $\omega(a)=\omega(\sigma(a))$.
\ebew

In the literature, a finite-dimensional algebra with a faithful functional is called a Frobenius algebra and the inverse of the automorphism $\sigma$ is called the Nakayama automorphism, see \cite{Na}. 
\ssnl
The terminology and notations we use come from the theory of locally compact groups and operator algebras. Then the automorphism $\sigma$ is called the modular automorphism and linear functionals admitting such a modular automorphism are called KMS-functionals, see e.g. \cite{Pe}.
\ssnl

We finish this item with a few remarks.

\opm
i) When the algebra $A$ has a unit, it is automatically non-degenerate and idempotent, i.e.\ $A^2=A$. 
\ssnl
ii) A faithful functional on an algebra $A$, finite-dimensional or not,  can only exist if the algebra is non-degenerate. Indeed, let $\omega$ be a faithful linear functional on $A$. Suppose that $a\in A$ and that $ac=0$ for all $c\in A$. Then $\omega(ac)=0$ for all $c$ and because $\omega$ is faithful, we have $a=0$. Similarly when $\omega(ca)=0$ for all $c$. Hence the product in $A$ is non-degenerate. 
\ssnl
iii) On the other hand, it will not be sufficient, even for a finite-dimensional algebra to have a unit when the product is non-degenerate as the  example below shows. In particular, there is no faithful functional in that case.
\eopm
In fact, with the following example, we see that we can have finite-dimensional \emph{idempotent} algebras, with a non-degenerate product, that are not unital.
\voorb
Let $A$ be the subspace of the algebra $M_n$ spanned by the matrix units $e_{1j}$ and $e_{jn}$ where $j$ runs from $1$ to $n$ in both cases. It is straightforward to verify that this is always an idempotent subalgebra of $M_n$, with a non-degenerate product. 
\ssnl
It is not unital when $n\geq 3$.
Consider e.g.\ the case with $n=3$.  If $x$ is $e_{11}, e_{13}, e_{23}$ or $e_{33}$ we have $xe_{23}=0$ while on the other hand $e_{12}e_{23}=e_{13}$. So $e_{23}$ does not belong to $Ae_{23}$. Therefore $A$ does not have a unit.
\evoorb
\nl
\bf Coproducts and counits \rm
\nl
We start with the following assumptions.

\voorw
Let $A$ be a finite-dimensional unital algebra. We assume that it has a coproduct. It is a liner map $\Delta:A\mapsto A\ot A$  
satisfying $(\Delta\ot\iota)\Delta=(\iota\ot\Delta)\Delta$ where $\iota$ is the identity map. We also need a counit. It is a linear map $\varepsilon:A\to \mathbb C$ satisfying $(\varepsilon\ot\iota)\Delta(a)=a$ and $(\iota\ot\varepsilon)\Delta(a)=a$ for all $a\in A$. 
\evoorw

We do not assume that the coproduct $\Delta$ is a homomorphism, only that it is coassociative. We also do not assume from the start that the counit is a homomorphism. Still, it is unique under the above conditions. 
\nl
Recall the following simple and well-known property.

\prop
Denote by $B$ the space of linear functionals on $A$ and use $\langle a,b\rangle$ to denote the value of $b$ in the point $a\in A$. The coproduct makes $B$ into an (associative) algebra by the formula 
$$\langle a,bb'\rangle=\langle \Delta(a),b\ot b'\rangle.$$ 
We use the obvious pairing of $A\ot A$ with $B\ot B$. The counit on $A$ is a unit for this algebra $B$ and 
$$\langle a,1\rangle=\varepsilon(a).$$
\eprop

Observe that the bilinear form is non-degenerate by definition. We use that for all $a\in A$ there is a linear functional $f$ on $A$ so that $f(a)\neq 0$.
\ssnl
Conversely, the product on $A$ induces a coproduct on $B$ by the formula 
$$\langle a\ot a',\Delta(b)\rangle=\langle aa',b\rangle.$$ The identity on $A$ gives a counit $\varepsilon$ on $B$ by the formula $\langle 1,b\rangle=\varepsilon(b)$.
\ssnl
From the point of view of duality, it is natural to assume that there is a counit for $(A,\Delta)$ as this is equivalent with the requirement that the dual algebra $B$ has a unit.
\ssnl
Recall what we mentioned already in the introduction, we use $\Delta$ both for the coproduct on $A$ and the coproduct on $B$ and similarly for the counit $\varepsilon$. 
\nl
\bf Dual pairs of algebras \rm
\nl
We see that an algebra $A$ with a coproduct $\Delta$ as above gives rise to a pairing of algebras as in the following definition.

\defin\label{defin:1.5}
Let $A$ and $B$ be algebras. A pairing of $A$ and $B$ is a non-degenerate bilinear form $(a,b)\mapsto \langle a,b\rangle$ from the Cartesian product $A\times B$ to $\mathbb C$.
\edefin

The coproduct  on one algebra is induced by the product on the other one. The counit  on one algebra is given by the  pairing with the identity of the other one. 
\ssnl
There is an obvious converse. Given such a pairing of algebras, the product in $B$ induces a coproduct on $A$ and the unit in $B$ gives a counit for this coproduct on $A$. Then the algebra $B$ is recovered as the dual of $A$, as described before.
\ssnl 
This observation is important for the rest of the paper and we explicitly summarize it in the following remark.

\opm
i) There are two obvious ways to deal with these objects. We can start with a pair of algebras $A$ and $B$ and then consider the coproducts and counits, coming from the pairing. We can also start with an algebra $A$, a coproduct $\Delta$ on this algebra and a counit $\varepsilon$ and consider the dual space $B$, endowed with the product and unit coming from the coproduct and counit of $A$. Clearly the two approaches are the same because we work with finite-dimensional algebras. 
\ssnl
ii) This is no longer true in infinite dimensions. Since our aim is to eventually extend the theory obtained here to infinite dimensions (cf. \cite{La-VD3b}), we mostly prefer to take a pair of algebras as the starting point. For some definitions, we take the other point of view. 
\eopm

Assume now that we have a dual pair of  algebras $A$ and $B$ as in Definition \ref{defin:1.5}. We consider the induced coproducts and counits as discussed before.

\prop\label{prop:1.7}
The counit on $B$ is a homomorphism if and only if $\Delta(1)=1\ot 1$ in $A\ot A$. Similarly for the counit on $A$.
\eprop

\bew
Take $b,b'\in B$. Assume that $\Delta(1)=1\ot 1$ in $A\ot A$. Then
\begin{equation*}
\varepsilon(bb')=\langle 1,bb'\rangle=\langle \Delta(1),b\ot b'\rangle=\langle 1,b\rangle \langle 1,b'\rangle
\end{equation*}
and we see that $\varepsilon(bb')=\varepsilon(b)\varepsilon(b')$. The converse is proven in the same way.
\ebew

For the further discussion, we do not assume that the coproducts are unital. But the condition will play a crucial role later.
\nl
We have \emph{associated actions} of one algebra on the other induced by the pairing as follows.

\prop\label{prop:1.8}
There are left and right actions of $A$ on $B$ and of $B$ on $A$ given by the formulas
\begin{align*}
\langle a,a'\tr b\rangle&=\langle aa',b\rangle = \langle a',b\tl a\rangle\\
\langle a\tl b,b'\rangle&=\langle a,bb'\rangle= \langle b'\tr a,b\rangle
\end{align*}
for all $a,a'\in A$ and $b,b'\in B$.
\eprop

There are some obvious properties of these actions. We have e.g. 
\begin{equation*}
\langle aa'a'',b\rangle=\langle aa',a''\tr b\rangle =\langle a,a'\tr(a''\tr b\rangle
\end{equation*}
so that $\langle a,(a'a'')\tr b\rangle=\langle a,a'\tr(a''\tr b)\rangle$. This holds for all $a$ and since the pairing is non-degenerate we find that $(a'a'')\tr b=a'\tr(a''\tr b)$. 

\prop\label{prop:1.9}
The four actions associated with a pairing of  algebras are faithful, unital and non-degenerate.
\eprop

\bew 
We prove this for the left action of $A$ on $B$. 
\ssnl
Suppose that $a\in A$ and that $a\tr b=0$ for all $b$. If we pair with the identity of $A$, we find that $\langle a,b\rangle=0$. This holds for all $b$ and because the pairing is non-degenerate, we have $a=0$. This shows that the action is faithful.
\ssnl
We clearly have $1\tr b=b$ for all $b$ so that the action is unital.
\ssnl
Finally, let $b\in B$ and assume that $a\tr b=0$ for all $a$. With $a=1$ we have $b=1\tr b=0$ and we see that the action is non-degenerate.
\ebew

This result is trivial but it will not be so in the general case of infinite-dimensional algebras in \cite{La-VD3b}. That is the reason for formulating it here as a separate result.

\nl
\bf Pairing of $^*$-algebras \rm
\nl
It is interesting to consider $^*$-algebras (i.e.\ algebras over $\mathbb C$ with an involution), not only because there are plenty of natural examples, but also because of the special role  of the antipode in this case. We have mentioned this already in the introduction, but now we give more details.
\ssnl
So assume now that $A$ and $B$ are $^*$-algebras and that we have a non-degenerate pairing of $A$ with $B$ as in Definition \ref{defin:1.5}.  

\prop\label{prop:1.10}
There is a bijective linear map $S_A$ from $A$ to $A$ satisfying 
$$\langle  S_A(a)^*,b\rangle=\overline{\langle a,b^*\rangle}$$
 for all $a,b$. Similarly, there is a bijective linear map $S_B$ from $B$ to itself satisfying 
 $$\langle a,S_B(b)^*\rangle=\overline{\langle  a^*,b\rangle}.$$ These two maps are each others adjoint, i.e.\ 
 $$\langle S_A(a),b\rangle=\langle a,S_B(b)\rangle$$ for all $a,b$.
\eprop

\bew
Let $a\in A$. The linear map $b\mapsto \overline{\langle a,b^*\rangle}$ is given by an element $a'$ in $A$. So $\langle a',b\rangle=\overline{\langle a,b^*\rangle}$ for all $b$. Define a linear map $S_A$ on $A$ by $S_A(a)={a'}^*$. Then we have 
\begin{equation*}
\langle  S_A(a)^*,b\rangle=\overline{\langle a,b^*\rangle}
\end{equation*}
for all $a\in A$ and $b\in B$. Similarly we have a linear map $S_B$ on $B$ defined by
\begin{equation*}
\langle  a,S_B(b)^* \rangle=\overline{\langle a^*,b\rangle}
\end{equation*}
for all $a,b$.
\ssnl
It is clear that the maps $a\mapsto S_A(a)^*$ and $b\mapsto S_B(b)^*$ are conjugate linear involutive maps. We also  find
\begin{equation*}
\langle S_A(a),b\rangle = \langle S_A(a)^{**},b\rangle=\overline{\langle S_A(a)^*,S_B(b)^*\rangle}
\end{equation*}
and by symmetry we get $\langle S_A(a),b\rangle=\langle a,S_B(b)\rangle$.
\ebew

In the operator algebraic approach there is a slightly different convention. The map $S_B$ is replaced by its inverse and so the last formula becomes $\langle S_A(a),b\rangle=\langle a,S_B\inv(b)\rangle$ instead, see e.g.\ \cite{VD-warsaw}.
We stick to the algebraic conventions as set out in the previous proposition.

\prop \label{prop:1.11}
We have $\Delta(S_A(a)^*)=\zeta((S_A\ot S_A)\Delta(a))^*$ for all $a$ where $\zeta$ is the flip map on $A\ot A$.  Similarly $\Delta(S_B(b)^*)=\zeta((S_B\ot S_B)\Delta(b))^*$. 
\eprop
\bew
Take $a\in A$ and $b,b'\in B$. Then
\begin{align*}
\langle \Delta(S_A(a)^*),b\ot b')
&=\langle S_A(a)^*,bb'\rangle \\
&=\langle a, (bb')^*\rangle^-\\
&=\langle a, {b'}^*b^* \rangle^- \\
&=\langle \Delta(a),  {b'}^*\ot b^* \rangle^- \\
&=\langle ((S_A\ot S_A)\Delta(a))^*,  b'\ot b \rangle. 
\end{align*}
We see that $\Delta(S_A(a)^*)=\zeta((S_A\ot S_A)\Delta(a))^*$. Similarly for $B$.
\ebew

 Recall that 
 $A\ot A$ and $B\ot B$ are $^*$-algebras with the natural involution defined by $(a\ot a')^*=a^*\ot {a'}^*$ and similarly for $B$. We also have used the short hand notation $\langle a,b\rangle^-$ for $\overline{\langle a,b\rangle}$.
 \ssnl
 Using a similar method as in the proof above, we get the following.

\prop\label{prop:1.12}
Assume that we have a pairing of two $^*$-algebras $A$ and $B$. The coproduct $\Delta$ on $A$ satisfies $\Delta(a^*)=\Delta(a)^*$ for all $a\in A$ if and only if  
$S_B(bb')^*=S_B(b)^*S_B(b')^*$ for all $b,b'\in B$. Similarly, for the coproduct on $B$
\eprop

\bew
Take $a\in A$ and $b,b'\in B$ and assume  that $S_B(bb')^*=S_B(b)^*S_B(b')^*$. Then we have
\begin{align*}
\langle a^*,bb'\rangle
&=\langle a,S_B(bb')^*\rangle^-\\
&=\langle a,S_B(b)^*S_B(b')^* \rangle^-\\
&=\langle \Delta(a),S_B(b)^*\ot S_B(b')^*\rangle^-\\
&=\langle \Delta(a)^*,b\ot b'\rangle. 
\end{align*}
For the last equality, we use the  formula for the involution on $A\ot A$. So we have $\Delta(a^*)=\Delta(a)^*$ for all $a$.
\ssnl
It is easy to see that the converse is also true. 
\ssnl
Similarly, $S_A(aa')^*=S_A(a)^*S_A(a')^*$ for all $a$ if and only if  $\Delta(b^*)=\Delta(b)^*$ for all $b$.

\ebew

If $\Delta$ is a $^*$-map on $B$, then $S_A$ is an anti-isomorphism of $A$. This follows from the above result. On the other hand, if $\Delta$ is a $^*$-map on $A$, then the result of Proposition \ref{prop:1.11} implies that $S_A$ also flips the coproduct. 
\ssnl
We are using the symbol $S$ because further it will behave like an antipode. That the maps $S_A$ and $S_B$ are anti-isomorphisms, flipping the coproduct is what we expect of an antipode.  However, the reader should be aware of the fact that this will not be enough to characterize the antipode. Indeed, if the algebra $A$ is abelian and $\Delta$ on $A$ is coabelian, also the identity map $\iota$ will have this property on $A$. 
\nl
These are  natural conditions for a pairing of $^*$-algebras, and hence it is reasonable to include them  as  conditions for a pairing of $^*$-algebras. This is done in the following definition.

\defin\label{defin:1.13}
Let $A$ and $B$ be $^*$-algebras. A non-degenerate bilinear form $(a,b)\mapsto \langle a,b\rangle$ on $A\times B$ is said to be a \emph{pairing of $^*$-algebras} if the coproducts are $^*$-maps.
\edefin

So, for such a pairing of $^*$-algebras, we have the linear maps  $S_A$ and $S_B$ on $A$ and $B$ resp.  They are anti-isomorphisms  and they flip the associated coproducts.
 \ssnl
 In this case, we moreover have the following property of the counits.

\prop
The counit $\varepsilon$ on $A$ satisfies $\varepsilon(a^*)=\overline{\varepsilon(a)}$ and $\varepsilon(S(a))=\varepsilon(a)$ for all $a$. Similarly for the counit on $B$.
\eprop

\bew
The first property follows from the uniqueness of the counit and because $\Delta$ is a $^*$-map. The second one is a consequence of the fact that $S_B$ is an anti-isomorphism of $B$ so that $S_B(1)=1$.
\ebew

\opm
If we start with a $^*$-algebra with a coproduct $\Delta$ satisfying $\Delta(a^*)=\Delta(a)^*$, we have the associated product on the dual $B$, but there is no canonical way to make $B$ into a $^*$-algebra. We need an anti-isomorphism $S_A$ on $A$ that flips the coproduct. Then we can define the involution on $B$ by the formula $\langle a,b^*\rangle=\langle S_A(a)^*,b\rangle^-$. We can define $S_B$ on $B$ as the adjoint of $S_A$ as before. Then we end up with a pairing of $^*$-algebras as in Definition \ref{defin:1.13}.
\eopm

Again in what follows, we will drop the indices and write $S$, both for $S_A$ on $A$ and for $S_B$ on $B$. 
\ssnl
We have the following formulas for the behavior of the maps $S$ with respect to actions, associated with the pairing of $^*$-algebras.

\prop\label{prop:1.16}
For all $a\in A$ and $b\in B$ we have 
\begin{equation*}
S(b\tr a)^*=S(a)^*\tl b^*
\tussenen
S(a\tr b)^*=S(b)^*\tl a^*.
\end{equation*}	 
\eprop

\bew
For another element $b'\in B$ we can write
\begin{align*}
\langle S(a)^*\tl b,b'\rangle 
&=\langle S(a)^*,bb')\rangle
=\langle a,(bb')^* \rangle^-\\
&=\langle a,{b'}^*b^*\rangle^-
=\langle b^*\tr a, {b'}^*\rangle^- \\
&=\langle S(b^*\tr a)^*,b'\rangle.
\end{align*}
We see that $S(a)^*\tl b=S(b^*\tr a)^*$. If we replace $b$ by $b^*$ we find $S(b\tr a)^*=S(a)^*\tl b^*$. 
\ssnl
The second property can be obtained  by a similar argument (or by symmetry).
\ebew

For these properties we use that $a\mapsto a^*$  and $b\mapsto b^*$ are involutions of the algebras. On the other hand, when we use that $S(aa')^*=S(a)^*S(a')^*$ and $S(bb')^*=S(b)^*S(b')^*$, we arrive at another pair of relations.

\prop\label{prop:1.17}
 For all $a\in A$ and $b\in B$ we have 
\begin{equation*}
 (a\tl b)^*=a^*\tl S(b)^*
 \tussenen
 (a\tr b)^*=S(a)^*\tr b^*.
\end{equation*}
\eprop
\bew
With $b'\in B$ we can write
\begin{align*}
\langle (a\tl b)^*,b' \rangle
&=\langle a\tl b,S(b')^*\rangle^-\\
&=\langle a,bS(b')^*\rangle^-\\
&=\langle a^*, S(b)^*b'\rangle\\
&=\langle a^*\tl S(b)^*,b'\rangle.
\end{align*}
We have used that $S(bb')^*=S(b)^*S(b')^*$ for all $b,b'$.
The proof of the second equality is obtained by using that $S(aa')^*=S(a)^*S(a')^*$ for all $a,a'$.
\ebew

The two sets of equalities are fundamentally different from each other. 
\ssnl 
Using that the antipodes are anti-isomorphisms of the algebras, we get similar forms for the action, not involving the involutions.

\prop\label{prop:1.18}
For all $a,b$ we have
\begin{equation*}
S(S(b)\tr a)=S(a)\tl b
\tussenen
S(S(a)\tr b)=S(b)\tl a.
\end{equation*}
\eprop

\bew
For $b'\in B$ we have
\begin{align*}
\langle S(S(b)\tr a),b'\rangle
&=\langle S(b)\tr a, S(b')\rangle \\
&=\langle a, S(b')S(b)\rangle \\
&=\langle a, S(bb')\rangle\\
&=\langle S(a),bb'\rangle \\
&=\langle S(a)\tl b, b'\rangle
\end{align*}
and so $S(S(b)\tr a)=S(a)\tl b$. The other equation is proved by using that $S$ is an anti-homomorphism on $A$.
\ebew

In fact, the result also follows by combining the formulas from the two previous results with the involutions. We have similar formulas with $S\inv$ because  $S\inv$ is also an anti-isomorphism. This results in the following analogous formulas
\begin{equation*}
a\tr b=S(S\inv(b)\tl S(a))
\tussenen
b\tr a=S(S\inv(a)\tl  S(b)).
\end{equation*}
\ssnl

These properties lead us naturally to the following section. But first, we introduce a new related concept.

\defin
Suppose that we have a pairing of algebras $A$ and $B$, not necessarily $^*$-algebras, as in Definition \ref{defin:1.5}. Assume that we have anti-isomorphisms $S_A$ on $A$ and $S_B$  on $B$, satisfying $\langle S_A(a),b\rangle=\langle a,S_B(b)\rangle$ for all $a,b$. Then we call $S_A$ and $S_B$ a pair of pre-antipodes. In this case, we say that the pair $(A,B)$ is \emph{a pair with pre-antipodes}.
\edefin

So, for a pairing of $^*$-algebras, as in Definition \ref{defin:1.13}, we always have pre-antipodes. 

\section{\hspace{-17pt}. Left and right integrals} \label{s:integrals} 

In this section we will introduce the notion of integrals. This is done in such a way that we get a nice and fairly general duality result. In the next section, we impose one more condition and get a duality of finite quantum hypergroups.
\ssnl
To explain where our definition of integrals comes from, we again look first at a pair of finite-dimensional unital $^*$-algebras as in Definition \ref{defin:1.13} of the previous section.
\nl
\bf Integrals  for a pairing of $^*$-algebras \rm
\nl
Consider  a non-degenerate pairing of finite-dimensional unital $^*$-algebras $A$ and $B$ as in Definition \ref{defin:1.13}. We have the anti-isomorphisms $S$ on $A$ and on $B$ as defined in Proposition \ref{prop:1.10}.
\ssnl
First recall the following. A linear functional $\varphi$ on a $^*$-algebra is called positive if $\varphi(a^*a)\geq 0$ for all $a$. A positive functional $\varphi$ is faithful if and only if $\varphi(a^*a)=0$ only holds for $a=0$. 
\ssnl
Assume that we have a faithful positive linear functional $\varphi$ on $A$. 
It gives a scalar product $(x,y)\mapsto \varphi(y^*x)$ where $x,y\in A$. 
Let $\mathcal H$ be the Hilbert space we get by completing $A$ for this scalar product and use $\Lambda$ for the injection of $A$ in $\mathcal H$. We denote the scalar product on $\mathcal H$ by $(\xi,\eta)\mapsto \langle \xi,\eta\rangle$. 
\ssnl
The following is standard.

\prop
Let $A$ act on $A$ by left multiplication.  Use $\pi(a)$ for the  linear map from $\Lambda(A)$ to itself given by $\pi(a)\Lambda(x)=\Lambda(ax)$. We clearly get a representation of $A$. In fact, we find that
\begin{equation*}
\langle \pi(a)\xi,\eta\rangle=\langle\xi,\pi(a^*)\eta \rangle
\end{equation*}
for $\xi$, $\eta$ in $\Lambda(A)$ and $a\in A$. So $\pi$ is a $^*$-representation. 
\eprop

In many cases, these operators $\pi(a)$ will be bounded and hence can be extended to operators on all of $\mathcal H$. In that case we have the GNS construction, see e.g.\ \cite{Pe}. However, in general, we can not expect to have a bounded representation. Fortunately, this is not important for our purposes.
\ssnl
On the other hand, consider the left action $(b,a)\mapsto a\tl S\inv(b)$ of $B$ on $A$. We use the pre-antipode $S$ on $B$ as obtained from the pairing (see Proposition \ref{prop:1.10}. This gives for each $b\in B$ an operator $\lambda(b)$ defined on $\Lambda(A)$ by $\lambda(b)\Lambda(a)=\Lambda(a\tl S\inv(b))$. 
Because the pre-antipode $S$ is assumed to be an anti-isomorphism, we get a representation of $B$. Also this operator will be  bounded  in many cases, but again, that is not important here.  Moreover, in many cases, this will be $^*$-representation. That however turns out to be important here 
due to the following result.

\prop\label{prop:2.2}
Let $b\in B$. Then $\langle\lambda(b)\xi,\eta\rangle=\langle\xi,\lambda(b^*)\eta\rangle$ for all $\xi,\eta$ in $\Lambda(A)$ if and only if
\begin{equation}
\varphi(c(a\tl S\inv(b)))=\varphi((c\tl b)a)\label{eqn:2.1}
\end{equation}
for all $a,c\in A$. 
\eprop

\bew
Let $a,c\in A$, $b\in B$ and assume that $\langle\lambda(b)\Lambda(a),\Lambda(c)\rangle=\langle\Lambda(a),\lambda(b^*)\Lambda(c)\rangle$. For the left hand side we find
\begin{equation*}
\langle \lambda(b)\Lambda(a),\Lambda(c)\rangle=\varphi(c^*(a\tl S\inv(b))).
\end{equation*}
For the right hand side we have
\begin{equation*}
\langle\Lambda(a),\lambda(b^*)\Lambda(c)\rangle
=\varphi((c\tl S\inv(b^*))^*a).
\end{equation*}
Now we use the first formula in Proposition \ref{prop:1.17}. We obtain 
$(c\tl S\inv(b^*))^*=c^*\tl b$.
If we insert this we find 
\begin{equation*}
\varphi(c^*(a\tl S\inv(b)))=\varphi((c^*\tl b)a).
\end{equation*}
Now we just have to replace $c$ by $c^*$ to get the required formula.
\ssnl
It is clear that the argument also works in the other direction. 
\ebew

If Equation \ref{eqn:2.1}  holds for all $b\in B$, then
\begin{equation}
(\iota\ot\varphi)((1\ot c)\Delta(a))=S((\iota\ot\varphi)(\Delta(c)(1\ot a)))\label{eqn:2.2}
\end{equation}
for all $a,c$.
Also conversely, if Equation (\ref{eqn:2.2}) holds for all $a,c$, then (\ref{eqn:2.1}) holds for all $a,b,c$.
\ssnl
This equality appears in many other theories of this kind. See e.g.\  (the proof of) Proposition 3.11 in \cite{VD-alg} and Proposition 1.5 in \cite{VD-W2}.
The above result shows why in the $^*$-algebra case this is a natural property. 
\nl
Equation (\ref{eqn:2.2}) no longer involves the  involutive structures.  Therefore it makes sense to impose it on the antipode in relation with the functional $\varphi$. It leads us to the notion of integrals in this framework.

\nl
\bf Left and right integrals \rm
\nl
Take a pairing of finite-dimensional unital algebras $A$ and $B$, not necessarily $^*$-algebras,  considered with the associate coproducts.
\ssnl
Motivated by the previous item, we introduce the main concept, namely left and right integrals. 

\defin\label{defin:2.3}
Let $\varphi$ be a linear functional on $A$. It is called \emph{a left integral} if there is a linear map
$S$ from $A$ to itself satisfying 
\begin{equation}
S((\iota\ot\varphi)(\Delta(a)(1\ot c)))=(\iota\ot\varphi)((1\ot a)\Delta(c))\label{eqn:2.3}
\end{equation} 
for all $a,c$ in $A$. Similarly a  linear functional $\psi$ is called \emph{a right integral} if there is a linear map $S'$  from $A$ to itself satisfying
\begin{equation}
S'((\psi\ot\iota)((c\ot 1)\Delta(a)))=(\psi\ot\iota)(\Delta(c)(a\ot 1))\label{eqn:2.4}
\end{equation}
for all $a,c$. 
\edefin

Remark that this definition is formulated for a pair $(A,\Delta)$ of an algebra $A$ with a coproduct $\Delta$ on $A$. We do not specify the map $S$. We will see below that, when $\varphi$ is faithful, $S$ is unique if it exists (see \ref{prop:2.5}).
\ssnl
We have the following equivalent formulation in terms of the pairing.

\prop\label{prop:2.4}
A linear functional $\varphi$ on $A$ is a left integral if there is a linear map $S$ on $B$ such that
\begin{equation}
\varphi((a\tl S(b))c)=\varphi(a(c\tl b))\label{eqn:2.5}
\end{equation}
for all $a,c\in A$ and $b\in B$. Similarly, a linear functional $\psi$ on $A$ is a right integral if there is a linear map $S'$ on $B$ such that
\begin{equation}
\psi(c(S'(b)\tr a))=\psi((b\tr c)a)
\end{equation}
for all $a,c\in A$ and $b\in B$.
\eprop

This is obtained from Equation (\ref{eqn:2.3} when pairing with an element $b$ of $B$. The map $S$ on $B$ is the adjoint of the map $S$ on $A$ that we have in the definition. Similarly for the map $S'$ on $B$.
\ssnl
We know that these properties are true for  integrals on Hopf algebras and multiplier Hopf algebras. There they are proven by using that the coproduct is a homomorphism. It is somewhat remarkable however that the results are still valid for quantum hypergroups, see \cite{De-VD1}, and also for weak Hopf algebras and weak multiplier Hopf algebras, see \cite{VD-W2}. When treating the examples in Sections \ref{s:hecke}, \ref{s:two-dim} and \ref{s:two-sub}, this will become more clear.
\opm
In the theory of operator algebras, there is an object called Kac algebras, see \cite{E-S}. It is an operator algebra with a coproduct. The existence of an antipode is assumed. It is an anti-isomorphism that flips the coproduct. In that theory, the left integral is defined with reference to the antipode using the same formula as above.
\eopm
An immediate and important consequence of the definition is the following.

\prop\label{prop:2.5}
For a faithful left integral $\varphi$, the associated map $S$ on $A$ is unique and injective. It flips the coproduct.  Similarly for the map $S'$ associated to a faithful right integral.
\eprop
\bew
i) We claim that elements of the form $(\iota\ot\varphi)(\Delta(a)(1\ot c))$ span all of $A$ when $\varphi$ is faithful. Indeed, suppose that $\omega$ is a linear functional on $A$ that kills all such elements. By the faithfulness of $\varphi$ it follows that 
$(\omega\ot\iota)\Delta(a)=0$ and if we apply the counit, we get $\omega(a)=0$ for all $a$ so that $\omega=0$. This proves the claim. 
\ssnl
We have a similar result for the three other types of elements we have in the Equations (\ref{eqn:2.3}) and (\ref{eqn:2.4}) in Definition \ref{defin:2.3}
\ssnl
ii) It follows that the map $S$ associated to $\varphi$ is uniquely defined and surjective. As the space is finite-dimensional, the map is bijective. 
\ssnl
iii) To prove that $S$ on $A$ flips the coproduct, we consider its adjoint on $B$ and show that this is an anti-isomorphism. For the adjoint $S$ on $B$ we have the formula
\begin{equation}
\varphi(a(c\tl b))=\varphi((a\tl S(b))c)\label{eqn:2.7}
\end{equation}
for all $a,c\in A$ and $b\in B$. Now take two elements $b,b'\in B$.
Then we have
\begin{align*}
\varphi((a\tl S(bb'))c)
&=\varphi(a(c\tl bb'))\\
&=\varphi((a\tl S(b'))(c\tl b))\\
&=\varphi((a\tl S(b')S(b))c).
\end{align*}
Because this holds for all $c$ and as $\varphi$ is faithful, we get $a\tl S(bb')=a\tl S(b')S(b)$. This holds for all $a$ and since the action is faithful, we find that $S(bb')=S(b)S(b')$.
\ssnl
iv) A similar argument works for the antipode $S'$ associated to a faithful right integral.
\ebew

Because of this result, the following definition makes sense.

\defin
Let $\varphi$ be a faithful left integral. We call $S$ the antipode associated with $\varphi$. Similarly we call $S'$ the antipode associated with a faithful right integral $\psi$.
\edefin

It is expected  that, under certain obvious conditions, a left integral, composed with its antipode, is a right integral with the same antipode. Indeed, and more precisely, we have the following result.

\prop\label{prop:2.7}
Assume that $\varphi$ is a left integral on $A$ with antipode $S$. Assume that $S$ is an anti-isomorphism of $A$.
Then $\varphi\circ S$ is a right integral on $A$ with the same antipode.
\eprop

\bew
This follows from the fact that $S$ flips the coproduct and that it is an anti-isomorphism. Indeed, define $\psi=\varphi\circ S$. Then, for all $a,c$ we have
\begin{align*}
S((\psi\ot\iota)((c\ot 1)\Delta(a)))
&=(\varphi\ot\iota)((S\ot S)(\Delta(a)(S(c)\ot 1)))\\
&=(\iota\ot\varphi)(\Delta(S(a))(\iota\ot S(c)))\\
&=S\inv((\iota\ot\varphi)(\iota\ot S(a))\Delta(S(c)))\\
&=S\inv((\varphi\ot\iota)((S\ot S)(\Delta(c)(a\ot 1))))\\
&=(\psi\ot\iota)(\Delta(c)(a\ot 1)).
\end{align*}

\vskip -20pt \ebew

We have a similar result for the composition $\psi\circ S'$ in the case of a right integral. The extra condition, namely that $S$ is an anti-isomorphism of $A$ is also natural as we have seen in the case of involutive algebras, see Proposition \ref{prop:1.12}.
\nl
\bf The duality theorem\rm
\nl
In this subsection, we construct the dual right integral on $B$ for a left integral on $A$.
\ssnl
Assume that we have a dual pair of algebras $A$ and $B$ as in Definition \ref{defin:1.5}. Assume that $\varphi$ is a faithful left integral on $A$ with antipode $S$. We consider also the adjoint of $S$ on $B$ and use again $S$ for this adjoint.

\lem\label{lem:2.8}
Let $b,d\in B$ and write $b=\varphi(\,\cdot\,c)$ for some $c\in A$. Then we have $db=\varphi(\,\cdot\,c')$ where $c'=c \tl S\inv(d)$.
\elem
\bew
For all $a$ in $A$ we have $\langle a,b\rangle=\varphi(ac)$ and by using the formula in Proposition \ref{prop:2.4}, that
\begin{equation*}
\langle a,db\rangle=\langle a\tl d,b\rangle = \varphi((a\tl d)c)=\varphi(a(c\tl S\inv(d))).
\end{equation*}
We see that $db=\varphi(\,\cdot\,c')$ where $c'=c\tl S\inv(d)$. 
\ebew

\prop\label{prop:2.9}
Define $\psi$ on $B$ by $\psi(b)=\varepsilon(c)$ if $b=\varphi(\,\cdot\,c)$ with $c\in A$. If the antipode $S$ associated with $\varphi$ is an anti-isomorphism of $A$, then $\psi$ is a faithful right integral on $B$ with antipode $S$.
\eprop

\bew
i) Because $\varphi$ is faithful, $\psi$ is well-defined on all of $B$.
\ssnl
ii) Let $b,d\in B$ and assume that $b=\varphi(\,\cdot\,c)$ for $c\in A$. From the lemma we see that $\psi(db)=\langle c,S\inv(d)\rangle$. Then the faithfulness of $\psi$ is a consequence of the fact that the pairing is non-degenerate.
\ssnl
iii) Let $a,c\in A$ and $b,d\in B$ with $b=\varphi(\,\cdot\,c)$. Then
\begin{equation*}
\psi(d(S(a)\tr b))=\langle S(a)c,S\inv(d)\rangle
\end{equation*}
because $S(a)\tr b=\varphi(\,\cdot\, S(a)c)$. On the other hand we have $a\tr d=S(S\inv(d)\tl S(a))$, see the remark after Proposition \ref{prop:1.18}. Then
\begin{align*}
\psi((a\tr d)b)
&=\langle c,S\inv(d)\tl S(a) \rangle\\
&=\langle S(a)c,S\inv(d)\rangle.
\end{align*}
We see that $\psi(d(S(a)\tr b))=\psi((a\tr d)b)$ so  $\psi$ is a right integral with antipode $S$.
\ebew

We again need that $S$ is an anti-isomorphism of $A$ to prove this result. That is natural because when we have an integral on $B$, we must have that its antipode is an anti-isomorphism on $A$, see Proposition \ref{prop:2.5}.
\ssnl
We find the following duality result.

\stel\label{stel:2.10}
Let $A$ and $B$ be unital finite-dimensional algebras with a pairing as in Definition \ref{defin:1.5}. Assume that there is a faithful left integral $\varphi$ on $A$ and that its antipode $S$ is an anti-isomorphism of $A$. Define $\psi$ on $B$ by $\psi(b)=\varepsilon(c)$ when $b=\varphi(\,\cdot\,c)$. Then $\psi$ is a faithful right integral on $B$ with antipode $S'$ on $B$, given as the adjoint of $S$ on $A$. The map $S'$ is an anti-isomorphism of $B$. 
\estel

If $b=\varphi(\,\cdot\,c)$, then $\psi(\,\cdot\,b)=S\inv(c)$, see Item ii) in the proof of Proposition \ref{prop:2.9}. It is an immediate consequence of the definition of $\psi$ as we see in Lemma \ref{lem:2.8}. Now we can take the dual left integral on $A$ constructed from $\psi$ on $B$ as above.  We get back the original left integral on $A$. We refer to this property as \emph{biduality}. 
\ssnl
The map $c\mapsto \varphi(\,\cdot\,c)$ can be considered as a generalization of the Fourier transform and  $b\mapsto \psi(\,\cdot\,b)$ (or rather its composition with $S$) as the inverse Fourier transform. We  refer to \cite{VD-ft} for a discussion on the Fourier transform in quantum group theory.
\ssnl
We also have the following alternative formulation of this duality result.

\prop\label{prop:2.12}
Let $(A,\Delta)$ be a pair of a finite-dimensional unital algebra $A$ with a coproduct $\Delta$. We assume that there is a counit on $A$ and that there is a faithful left integral on $A$ with an antipode that is an anti-isomorphism. Then the dual pair $(B,\Delta)$ is again a unital algebra with a coproduct, admitting a counit. Moreover there is a left integral on $B$ with an antipode that is an anti-isomorphism.
\eprop

To show this, we use that a right integral, composed with its antipode is a left integral when this antipode is an anti-isomorphism and we apply the theorem.
\ssnl
Recall that the antipodes will flip the coproducts by the very existence of the integrals.
\ssnl
We call $(B,\Delta)$ the dual of $(A,\Delta)$. The dual of $(B,\Delta)$ is the same as the original  pair $(A,\Delta)$.
\ssnl
In the case of a $^*$-algebra we have the following result.

\prop\label{prop:2.13}
If $\varphi$ is a left integral, $a\in A$ and $b=\varphi(\,\cdot\,a)$, then $\psi(b^*b)=\overline{\varphi(a^*a)}$. In particular, if $\varphi$ is positive, then so is $\psi$ and $\psi(b^*b)=\varphi(a^*a)$.
\eprop
\bew
We know that $\psi(b^*b)=\langle S\inv(a),b^*\rangle$ and so
\begin{equation*}
\psi(b^*b)=\langle a^*,b\rangle^-=\varphi(a^*a)^-.
\end{equation*}
The second statement is obvious.
\ebew

If we think of the map $a\mapsto \varphi(\,\cdot\,a)$ as a Fourier transform, then the above result is like Plancherel's theorem.  Again see \ \cite{VD-ft}.
\ssnl
 We have seen that the composition $\varphi\circ S$ is a right integral if $S$ is an anti-isomorphism. If $S$ is a $^*$-map, then this right integral will again be positive when $\varphi$ is positive. However, it is not expected that this is true in general.

\section{\hspace{-17pt}. Finite quantum hypergroups} \label{s:finite}

In what follows, we again take finite-dimensional unital algebras $A$ and $B$, not necessarily $^*$-algebras, with a non-degenerate pairing $(a,b)\mapsto \langle a,b\rangle$ from $A\times B$ to $\mathbb C$ as in Definition \ref{defin:1.5}. In the previous section, we obtained a duality of algebras with integrals under the extra condition that the associated antipodes are anti-isomorphisms. The intention here is to impose more conditions on the pairing to get a pair of \emph{finite quantum hypergroups}.
\nl
\bf Left and right invariant functionals \rm
\nl
We begin with the introduction of  \emph{left} and \emph{right invariant functionals} for a pair $(A,\Delta)$ of a finite-dimensional unital algebra with a coproduct. Recall that we do not assume that $\Delta$ is a homomorphism. The condition $\Delta(1)=1\ot 1$ is also not required from the beginning.

\defin\label{defin:3.1}
A linear functional $\varphi$ on $A$ is called \emph{left invariant} if $(\iota\ot \varphi)\Delta(a)=\varphi(a)1$ for all $a\in A$. Similarly, a linear functional $\psi$ is called \emph{right invariant} if $(\psi\ot\iota)\Delta(a)=\psi(a)1$ for all $a$.
\edefin

In the previous section, we have introduced the notion of 
a left integral and of a right integral.   There are examples where invariant functionals exist that are not integrals. We give such examples at the end of Section \ref{s:two-dim}.  See also some comments in the final section \ref{s:conclusions} where we  draw conclusions. 
\ssnl
Because the algebras are finite-dimensional, a linear functional $\varphi$ on $A$ is given by an element $h\in B$ satisfying $\varphi(a)=\langle a,h\rangle$ for all $a$. Left invariance of $\varphi$ is equivalent with the property that $bh=\varepsilon(b)h$ for all $b\in B$. In Hopf algebra theory, an element with this property is called a cointegral (or an integral \emph{in} the algebra). 
\ssnl
As it is our intention to generalize our results to the infinite-dimensional case (see\cite{La-VD3b}), we will work with the integrals and not with the cointegrals. The existence of such an element in $B$ is only guaranteed in finite dimensions. 
\ssnl

The following result gives  a necessary condition to have invariant functionals. 

\prop\label{prop:3.2}
If we have a non-zero  left invariant linear functional on $(A,\Delta)$, then $\Delta(1)=1\ot 1$.
\eprop

\bew
Let $\varphi$ be a left invariant functional. Take any $a\in A$ and assume that $\varphi(a)\neq 0$. Then
\begin{align*}
\Delta(1)\varphi(a)
&=\Delta((\iota\ot\varphi)(\Delta(a)))\\
&=(\iota\ot\iota\ot\varphi)((\Delta\ot\iota)\Delta(a))\\
&=(\iota\ot\iota\ot\varphi)((\iota\ot\Delta)\Delta(a))\\
&=(\iota\ot\iota\ot\varphi)\Delta_{13}(a)=\varphi(a)\,1\ot 1.
\end{align*}
We use the \emph{leg numbering notation}. More precisely $\Delta_{13}(a)=(\iota\ot\zeta)(\Delta(a)\ot 1)$  where $\zeta$ is the flip map on $A\ot A$. 
This proves the result.
\ebew

We see that a non-trivial left invariant functional can only exist when $\Delta$ is unital. Hence, in the case of a finite quantum groupoid, we have integrals that are not invariant functionals. See Example \ref{voorb:5.11}.
\ssnl
Recall from Proposition \ref{prop:1.7} that the coproduct on $A$ is unital if and only if the counit $\varepsilon$ on $B$ is a homomorphism. Therefore, we can only have invariant functionals on $A$ when the counit on $B$ is a homomorphism. 
In fact, there is the following inverse. 

\prop\label{prop:3.3}
Assume that the counit on $B$ is a homomorphism.
If there is a \emph{faithful} linear functional on $B$, then there is a left and a right invariant functional on $A$.
\eprop

\bew
Let $\omega$ be a faithful linear functional on $B$. Because $B$ is finite-dimensional, by Proposition \ref{prop:1.2}, there is an element $h\in B$ so that $\varepsilon=\omega(\,\cdot\,h)$ for the counit $\varepsilon$ on $B$. 
Then, for all $b,b'\in B$, we have 
\begin{equation*}
\omega(b'bh)=\varepsilon(b'b)=\varepsilon(b')\varepsilon(b)=\varepsilon(b)\omega(b'h). 
\end{equation*}
This holds for all $b'$. Because $\omega$ is faithful, we have $bh=\varepsilon(b)h$ for all $b$. This means that $\varphi$, defined on $A$ by $\varphi(a)=\langle a,h\rangle$, is left invariant. 
\ssnl
Similarly, if we take $k\in B$ satisfying $\varepsilon=\omega(k\,\cdot\,)$, we will have $kb=\varepsilon(b)k$ for all $b\in B$. Then $a\mapsto \langle a,k\rangle$ is a right integral.
\ebew

When $B$ is an operator algebra,  i.e.\ a $^*$-algebra with a faithful \emph{positive} linear functional, we can improve this result.

\prop\label{prop:3.4}
Again assume that the counit on $B$ is a homomorphism.
Now also assume that $B$ is an operator algebra. Then there is a self-adjoint idempotent element $h\in B$ satisfying  $bh=hb=\varepsilon(b)h$ for all $b\in B$. It defines a linear functional $\varphi$ on $A$ by $\varphi(a)=\langle a,h\rangle$. It is left and right invariant and has the property that $\varphi(1)=1$. It is the unique invariant functional on $A$ with $\varphi(1)=1$.
\eprop

\bew
i) Because the counit is a homomorphism, its kernel is a two-sided ideal of $B$. As $B$ is an operator algebras it is a direct sum of matrix algebras. Hence there is a self-adjoint idempotent $h$ with the property that $bh=hb=0$ when $\varepsilon(b)=0$. This implies that $bh=hb=\varepsilon(b)h$ for all $b\in B$. We also have $\varepsilon(h)=1$. 
\ssnl
ii) The element $h$ gives a functional $\varphi$ on $A$ that is both left and right invariant. Moreover $\varphi(1)=1$.
\ssnl
iii) Finally, for any left invariant functional $\varphi'$ on $A$ we have that for all $a$
\begin{equation*}
(\varphi\ot\varphi')\Delta(a)=\varphi'(a)
\tussenen
(\varphi\ot\varphi')\Delta(a)=\varphi'(1)\varphi(a).
\end{equation*}
For the first equality we use that $\varphi'$ is  left invariant and  that $\varphi(1)=1$. For the second one we use that $\varphi$ is right invariant. This shows that in this case invariant functionals are unique (up to a scalar).
\ebew

\ssnl
It is not possible to conclude that the invariant functional on $A$ we get in this way is faithful. 

\opm
If the coproduct is a homomorphism, we can construct an antipode if we have faithful invariant functionals. This follows from the Larson Sweedler theorem. For the case of a unital coproduct,  see \cite{La-Sw}. For the more general case, see \cite{VD-W1}. 
In this case, a left invariant functional is also a left integral. In particular, we get uniqueness of integrals. 
\eopm

As we aim to study quantum hypergroups, we do not require that the coproduct is a homomorphism. Then, unfortunately, we can not say much more at this point. We see from Proposition \ref{prop:3.3} that invariant functionals exist under mild conditions. To have integrals, we need the existence of an antipode as in Definition \ref{defin:2.3}. 
\nl
\bf Invariant functionals and integrals \rm
\nl
Recall that there is a difference between invariant functionals and integrals. Here we study the relation between the two concepts.
\ssnl
The following is easy to prove.

\prop\label{prop:3.5}
Assume that $(A,\Delta)$ admits a left integral $\varphi$. Denote by $S$ the associated antipode. If $\Delta$ is unital and if $S(1)=1$, then $\varphi$ is left invariant. Similarly for a right integral.
\eprop

\bew
With $c=1$ Equation (\ref{eqn:2.3}) in Definition \ref{defin:2.3} reads as 
\begin{equation*}
 S((\iota\ot\varphi)\Delta(a))=(\iota\ot\varphi)((1\ot a)\Delta(1))=\varphi(a)1.
\end{equation*}
 Then the result follows because $S(1)=1$ so that also $S\inv(1)=1$. Similarly for a right integral.
\ebew

That $\Delta$ has to be unital is obviously needed for this property as non-trivial invariant functionals can only exits in that case, see Proposition \ref{prop:3.2}.

\opm
Suppose that we do not  assume that $S(1)=1$ in the above proposition. Then we can apply 
the basic equality in Definition \ref{defin:2.3} with $c=S(1)$. Because  $\Delta$ is unital and $S$ flips the coproduct, we get $\Delta(S(1))=S(1)\ot S(1)$. Then
\begin{equation*}
S((\iota\ot\varphi)(\Delta(a)(1\ot S(1))))=(\iota\ot\varphi)((1\ot a)\Delta(S(1)))=\varphi(aS(1))S(1)
\end{equation*}
for all $a$ and we see that the functional $a\mapsto \varphi(aS(1))$ is left invariant. This is a slight improvement of the above proposition.
\eopm

The condition $S(1)=1$ will follow  when we require that $S$ is an anti-isomorphism of the algebra $A$. We know that also this is a natural condition. We needed it to have the duality in the previous section and we will in the end also need it here.
\ssnl
But first we prove the following \emph{uniqueness property}.

\prop\label{prop:3.7}
Assume that $\varphi$ is a faithful left integral on $A$ with antipode $S$ and that $S(1)=1$. 
Let $\varphi'$ be any left invariant functional on $A$. Then $\varphi'$ is a scalar multiple of $\varphi$.
\eprop

\bew
By Proposition \ref{prop:1.2} we have an element $c$ so that $\varphi'(a)=\varphi(ac)$ for all $a$. Then we have
\begin{align*}
\varphi(ac)1
&=\varphi'(a)1\\
&=(\iota\ot \varphi')\Delta(a)\\
&=(\iota\ot\varphi)(\Delta(a)(1\ot c))\\
&=S\inv((\iota\ot\varphi)((1\ot a)\Delta(c)))
\end{align*}
for all $a$.  By the faithfulness of $\varphi$ we get $S(1)\ot c=\Delta(c)$. Now apply the counit and we get $c=\varepsilon(c)S(1)$. 
Since $S(1)=1$, the result follows.
\ebew

We cannot get uniqueness using the argument in the proof of Proposition \ref{prop:3.4} because it can happen that $\varphi(1)=0$.
\ssnl
Without the requirement that $S(1)=1$, we would get that $\varphi'(a)=\lambda\varphi(aS(1))$ with $\lambda=\varepsilon(c)$. This gives uniqueness of invariant functionals if they exist when there is a faithful left integral. Observe that we have shown earlier that $a\mapsto \varphi(aS(1))$ is left invariant when it is given that $\Delta$ is unital. These two results are therefore in agreement with each other.
\ssnl
The uniqueness result proven in Proposition \ref{prop:3.7} has some important consequences, see further. We first give another proof, providing also some extra results.

\prop\label{prop:3.8}
Assume that we have a faithful left integral $\varphi$ on $A$ with an antipode $S$ satisfying $S(1)=1$. Assume that $\varphi'$ is a left invariant functional. Then there is an element $\delta$ in $A$ so that $(\varphi'\ot\iota)\Delta(a)=\varphi(a)\delta$. We also have $\varphi'(S(a))=\varphi(a\delta)$ for all $a$.
\eprop

\bew
Assume that $\varphi$ is a left integral and that $\varphi'$ satisfies $(\iota\ot\varphi')\Delta(a)=\varphi'(a)1$. 
As $\varphi$ is a left integral with antipode $S$ we get
\begin{equation*}
(\varphi'\ot \varphi)((1\ot a)\Delta(c))=(\varphi'\circ S\ot\varphi)\Delta(a)(1\ot c)).
\end{equation*}
Because $S$ flips the coproduct and $S(1)=1$ we have that $\varphi'\circ S$ is right invariant. Then for the right hand side we find $\varphi'(S(a))\varphi(c)$. For the left hand side we find $\varphi(a\delta_c)$ where $\delta_c=(\varphi'\ot\iota)\Delta(c)$.  Hence $\varphi'(S(a))\varphi(c)=\varphi(a\delta_c)$.
\ssnl
Now choose $c_0$ so that $\varphi(c_0)=1$ and use $\delta$ for the associated element. Then $\varphi'(S(a))=\varphi(a\delta)$.  If we insert this in the original formula, we find $\varphi(c)\varphi(a\delta)=\varphi(a\delta_c)$. This is true for all $a$ and because $\varphi$ is faithful we get
$\delta_c=\varphi(c)\delta$ for all $c$. This means that $(\varphi'\ot\iota)\Delta(c)=\varphi(c)\delta$.
\ebew

If we apply the counit on the last equation, we find $\varphi'(c)=\varphi(c)\varepsilon(\delta)$.
\ssnl
We can apply this with $\varphi$ in the place of $\varphi'$ if we include the condition that $\Delta(1)=1\ot 1$. Indeed, then we know by Proposition \ref{prop:3.5} that $\varphi$ is left invariant. We obtain the existence of an element $\delta$ in $A$ satsfying $(\varphi\ot\iota)\Delta(a)=\varphi(a)\delta$ and $\varphi(S(a))=\varphi(a\delta)$. This property is also true for algebraic quantum groups, see \cite{VD-alg} and for algebraic quantum hypergroups, see \cite{De-VD1}. It is called the modular element. In fact, the proof above is inspired by the arguments used to obtain these results.

\prop
Assume that $\Delta$ is unital and that  there is a faithful left integral $\varphi$ with an antipode $S$ satisfying $S(1)=1$. Then any left invariant functional is a scalar multiple of $\varphi$. In particular, any non-zero left invariant functional is a faithful left integral. Moreover, any non-zero left integral is a faithful left invariant functional.
\eprop

\bew
We have seen in Proposition \ref{prop:3.7} that any left invariant functional is a scalar multiple of $\varphi$. Consequently, if it is non-zero, it is itself a faithful left integral. On the other hand, any non-zero left integral is left invariant by Proposition \ref{prop:3.5} and hence by the previous result, it has to be faithful.
\ebew

So we get uniqueness of left integrals if $\Delta$ is unital and if there is a left integral with an antipode satisfying $S(1)=1$. Because the antipode is determined by the integral, we also have a unique antipode for the pair $(A,\Delta)$. 

\nl
\bf Finite quantum hypergroups \rm
\nl
If we want duality for finite quantum hypergroups, we will need to impose  the same conditions on $A$ that we have on $B$. In the first place, we need a faithful integral on both so that we have a pair like in Theorem \ref{stel:2.10}. In particular, we need antipodes that are anti-isomorphisms. 
Moreover, and this \emph{distinguishes the finite quantum hypergroups} from the objects we had in the previous section, we need unital coproducts on both sides. Or equivalently, we need that the counits are homomorphisms. 
\ssnl
This leads to the following definition.

\defin\label{defin:3.10}
Let $A$ be a finite-dimensional unital algebra with a coproduct $\Delta$ that admits a counit. Assume that there is a faithful left integral on $A$. If the coproduct is unital, the counit a homomorphism and the antipode associated with the integral an anti-isomorphism, then we call $(A,\Delta)$ a  \emph{finite quantum hypergroup}.
When $A$ is a $^*$-algebra and $\Delta$ a $^*$-map, we call $A$ a \emph{finite $^*$-quantum hypergroup}.
\edefin

Remember that the existence of a faithful left integral implies that its antipode $S$ also flips the coproduct. 
\ssnl
Compare this with Definition 1.10 from \cite{De-VD1}. There the invariance of $\varphi$ is defined by the equality $(\iota\ot\varphi)(\Delta(a)=\varphi(a)1$ and the antipode is defined as we do here. It follows that $\Delta$ must be unital. And we see that we precisely obtain the same notion for a finite-dimensional algebra.
\ssnl
Note that the counit is unique if it exists.  As we assume that the antipode $S$ of the left integral is an anti-isomorphism, we must have $S(1)=1$.  Then we also have a unique left integral and the antipode associated with it is also unique. Therefore we do not need to include these objects in the notation.
\nl
From the duality results in the previous section, when applied to a finite quantum hypergroup, we find the following.

\stel\label{stel:3.10}
Let $(A,\Delta)$ be a finite quantum hypergroup. Consider the dual space $B$ with product and coproduct $\Delta$ adjoint to the coproduct and product of $A$. Then $(B,\Delta)$ is again a finite quantum hypergroup. 
\estel

\bew
i) The algebra $B$ is unital. We have a coproduct $\Delta$ on $B$. It admits a counit, given by the identity in $A$. The coproduct on $B$ is unital because the counit on $A$ is assumed to be an homomorphism. On the other hand, the counit on $B$ is a homomorphism because the coproduct on $A$ is unital.
\ssnl
ii) Let $\varphi$ be a faithful left integral on $A$ and $S$ its antipode. By Proposition \ref{prop:2.9} we have a faithful right integral $\psi$ on $B$. Its antipode is the adjoint of the antipode on $A$. It is an anti-isomorphism of $B$ because the antipode on $A$ flips the coproduct (by Proposition \ref{prop:2.5}). The antipode on $B$ flips the coproduct as the antipode on $A$ is an anti-isomorphism.
\ssnl
iii) If we compose $\psi$ with its antipode, we get a faithful left integral on $B$. Hence $(B,\Delta)$ is again a finite-quantum hypergroup.
\ebew

The biduality result that we have for finite-dimensional unital algebras with a faithful integral (see Theorem \ref{stel:2.10} and Proposition \ref{prop:2.12}) is also valid here. The extra assumption, namely that the coproducts are unital, is satisfied for the algebra and its dual.
\ssnl
We end this section with a brief look at the involutive case. We have the following version of Theorem \ref{stel:3.10} for $^*$-algebras.

\stel
Assume that $(A,\Delta)$ is a finite $^*$-quantum hypergroup with antipode $S$. The dual is a $^*$-algebra if the involution is defined by $\langle a,b^*\rangle=\langle S(a)^*,b\rangle^-$. It is again a $^*$-quantum hypergroup. If the left integral on $A$ is positive, then so is the dual right integral on $B$.
\estel

For the last statement, see Proposition \ref{prop:2.13}.
\ssnl
In the next three sections, we will treat several examples to illustrate the notions and results obtained in the preceding  sections.

\section{\hspace{-17pt}. Classical Hecke algebras} \label{s:hecke}

In the first place, we look at the case of a finite group with a subgroup. This gives rise to the classical Hecke algebras.  But we go beyond these and derive some other, related examples in the next section.
\ssnl
Most calculations are easy and therefore we leave them as an exercise for the reader. 
\nl
\bf The basic examples associated with a group $G$ and a subgroup $H$ \rm
\nl
Let $G$ be a finite group. First we recall the following well-known fact.
\ssnl
Let $A_0$ be the $^*$-algebra of all complex functions $f$ on $G$ with pointwise operations. It is a Hopf $^*$-algebra for the coproduct $\Delta_0$ on $A_0$ defined by $\Delta_0(f)(p,q)=f(pq)$ where $p,q\in G$. The counit is given by $\varepsilon(f)=f(e)$ and the antipode by $(S(f))(p)=f(p\inv)$ for all $p$. 
\ssnl
On the other hand, let $B_0$ be the group algebra $\mathbb C G$. It is also a Hopf $^*$-algebra. If we use $p\mapsto \lambda_p$ for the canonical embedding of $G$ in the group algebra, we have $\lambda_p\lambda_q=\lambda_{pq}$ and 
$\lambda_p^*=\lambda_{p\inv}$ for all $p\in G$. The coproduct $\Delta_0$ on $B_0$ is given by $\Delta_0(\lambda_p)=\lambda_p\ot\lambda_p$ and the counit satisfies $\varepsilon(\lambda_p)=1$ for all $p$. For the antipode we have $S(\lambda_p)=\lambda_{p\inv}$.
\ssnl
There is a natural pairing between $A_0$ and $B_0$ given by $\langle f, \lambda_p\rangle=f(p)$ for all $p\in G$. It is a pairing of Hopf $^*$-algebras in the sense that $B_0$ is the dual of $A_0$.
\nl
Now let $H$ be a subgroup of $G$. We will associate a dual pair of finite $^*$-quantum hypergroups $A$ and $B$. They are derived from the Hopf $^*$-algebras $A_0$ and $B_0$ using a similar procedure.

\defin 
We define the linear map $E$ from $A_0$ to itself by 
\begin{equation*}
(E(f))(p)=\frac1{n^2}\sum_{h,k\in H}f(hpk)
\end{equation*}
for $p\in G$. Here $n$ is the number of elements in $H$. We use $A$ for the range of $E$.
\edefin

\prop\label{prop:4.2} 
The subspace $A$ is the $^*$-subalgebra of $A_0$ of functions $f$ constant on double cosets. Moreover $E$ is a conditional expectation from $A_0$ to $A$, it is unital and self-adjoint.
\eprop

\bew
This is mostly trivial. That $E$ is a conditional expectation means that $E(f)=f$ when $f\in A$ and that $E(ff')=E(f)f'$ and $E(f'f)=f'E(f)$ when $f\in A_0$ and $f'\in A$.
\ebew

It is clear that $E(S(f))=S(E(f))$ for all $f$. On the other hand, for the counit we have $\varepsilon(E(f))=\frac1{n}\sum_{h\in H}f(h)$ and this is not $\varepsilon(f)$ (except when $H$ is the trivial subgroup $\{e\}$).  
\ssnl
For the coproduct, we have the following crucial property.

\prop
If $f\in A$ then 
\begin{equation*}
(E\ot E)\Delta_0(f)(p,q)=(E\ot\iota)\Delta_0(f)(p,q)=(\iota\ot E)\Delta_0(f)(p,q)=\frac1{n}\sum_{k\in H}f(pkq).
\end{equation*}
These elements belong to $A\ot A$. 
\eprop

\bew
For $f\in A$ we have, for all $p,q\in G$,
\begin{align*}
((E\ot\iota)\Delta_0(f))(p,q)
&=\frac1{n^2}\sum_{h,k\in H}\Delta_0(f)(hpk,q)\\
&=\frac1{n^2}\sum_{h,k\in H}f(hpkq)\\
&=\frac1{n}\sum_{k\in H} f(pkq).
\end{align*}
Similarly $((\iota\ot E)\Delta_0(f))(p,q)=\frac1{n}\sum_{k\in H}f(pkq)$. We see that $(E\ot\iota)\Delta_0(f)=(\iota\ot E)\Delta_0(f)$. Hence these functions belong to $A\ot A$. So we also have 
$$(E\ot E)\Delta_0(f)=(E\ot\iota)\Delta_0(f)=(\iota\ot E)\Delta_0(f)=\frac1{n}\sum_{k\in H}f(pkq).$$
\ebew

We now have the main property.

\stel\label{stel:4.4}
Define $\Delta$ on $A$ by $\Delta(f)=(E\ot \iota)\Delta_0(f)$. Then $(A,\Delta)$ is a finite $^*$-quantum hypergroup.
\estel

\bew
i) First we argue that the coproduct $\Delta$ on $A$ is coassociative. Indeed, for all $f\in A$ we have
\begin{align*}
(\Delta\ot\iota)\Delta(f)&=(E\ot \iota \ot E)((\Delta_0\ot\iota)\Delta_0(f))\\
(\iota\ot\Delta)\Delta(f)&=(E\ot \iota \ot E)((\iota\ot\Delta_0)\Delta_0(f))
\end{align*}
and it follows from coassociativity of $\Delta_0$. We also have $\Delta(1)=1\ot 1$ because $E(1)=1$. Further $\Delta$ is a $^*$-map on $A$ because that is so for $\Delta_0$ and $E$ is self-adjoint.
\ssnl
ii) For the counit $\varepsilon$ on $A$ we find, when $f\in A$, that 
\begin{equation*}
(\varepsilon\ot\iota)\Delta(f)=E((\varepsilon\ot\iota)\Delta_0(f))=E(f)=f
\end{equation*}
and similarly for the other side. So the restriction of the counit to $A$ is a counit on $A$ for $\Delta$. It is a homomorphism.
\ssnl
iii) Define $\varphi(f)=\sum_{p\in G} f(p)$. Let $S$ be the antipode of the Hopf algebra $(A_0,\Delta_0)$. Then we have
\begin{equation*}
S((\iota\ot\varphi)(\Delta_0(f)(1\ot f')))=(\iota\ot\varphi)((1\ot f)\Delta_0(f'))
\end{equation*}
for all $f,f'\in A_0$. This is a well-known property of the antipode and in this case easy to verify. Now suppose that $f,f'\in A$. Then apply $E$ on the first factor in this equation and use that $E\circ S=S\circ E$. We find
\begin{equation*}
S((\iota\ot\varphi)(\Delta(f)(1\ot f')))=(\iota\ot\varphi)((1\ot f)\Delta(f')).
\end{equation*}
This proves that the restriction of $\varphi$ to $A$ is a left integral with $S$ on $A$ as its antipode. 
\ssnl
iv) The restriction of $\varphi$ to $A$ is still faithful. In this case we can use that $\varphi$ is positive and faithful on $A_0$ and therefore its restriction to the $^*$-subalgebra $A_0$ is still faithful.
\ssnl 
v) The antipode is an anti-homomorphism on $A_0$ and so also on $A$. It flips the coproduct $\Delta_0$ and hence
\begin{equation*}
\Delta(S(f))=(E\ot\iota)\zeta(S\ot S)\Delta_0(f)= \zeta(S\ot S)(\iota\ot E)\Delta_0(f)=\zeta(S\ot S)\Delta(f)
\end{equation*}
because $E$ and $S$ commute. We use $\zeta$ for the flip map on $A_0\ot A_0$. So $S$ also flips the coproduct $\Delta$ on $A$.
\ebew

It is also possible to give a proof of the main property of the antipode on $A$ without using the connection with $\Delta_0$. For this, take $f,g\in A$ and $p\in G$. Then 
\begin{align*}
((\iota\ot\varphi)(\Delta(f)(1\ot g)))(p)
&=\sum_q \Delta(f)(p,q)g(q) \\
&=\frac{1}{n}\sum_{h,q} f(phq)g(q) \\
&=\frac{1}{n}\sum_{h,q} f(q)g(h\inv p\inv q) \\
&=\sum_q f(q)g(p\inv q).
\end{align*}
We take the sum over $h\in H$ and $q\in G$.
On the other hand 
\begin{equation*}
((\iota\ot\varphi)(1\ot g)\Delta(f))(p\inv)=\frac{1}{n}\sum_{h,q} g(q)f(p\inv hq)=\sum_q g(q)f(p\inv q).
\end{equation*}
This shows that 
\begin{equation*}
S((\iota\ot\varphi)(\Delta(f)(1\ot g)))=(\iota\ot\varphi)((1\ot g)\Delta(f))
\end{equation*}
for all $f,g\in A$ and therefore $S$ is an antipode for $\varphi$.
\ssnl

\opm
By using the conditional expectation $E$, we can use the knowledge about the original Hopf algebra. Moreover, it turns out to be a general procedure to get a quantum hypergroup from a quantum group, see e.g. \cite{De-VD2}.  In particular, we see why the basic formula relating the integral with the antipode still holds. We will also use it for the dual quantum hypergroup of this example. Finally, the use of a conditional expectation to construct quantum hypergroups will also be illustrated in  Section \ref{s:two-sub}. See e.g.\ Theorem \ref{stel:6.16}.
\eopm

When $H$ is a normal subgroup of $G$ we get a Hopf $^*$-algebra, but as we see from the next result, only in that case.

\prop\label{prop:4.6a}
The coproduct  is a homomorphism on $A$ if and only if $H$ is a normal subgroup of $G$.
\eprop

\bew
i) Assume that $H$ is a normal subgroup. Then we have, for $f\in A$, 
\begin{equation*}
\Delta(f)(p,q)=\frac1{n}\sum_{h\in H}f(phq)=\frac1{n}\sum_{h\in H}f(php\inv pq)=f(pq)
\end{equation*}
because $php\inv\in H$ for all $p\in G$ and $h\in H$. Therefore in this case we have $\Delta=\Delta_0$ on $A$. Hence $\Delta$ is a homomorphism on $A$. 
\ssnl
ii) Conversely assume that $\Delta$ is a homomorphism. Take a point $p\in G$ and consider the map $f\mapsto \Delta(f)(p,p\inv)$. It is a homomorphism from $A$ to $\mathbb C$ and so there is an element $q$ in $G$ such that $\Delta(f)(p,p\inv)=f(q)$ for all $f\in A$. 
We use that $A$ is the function algebra on the set of double cosets of $G$ over $H$. It is an abelian algebra and any homomorphism  on it is given by evaluation in a point of the set. Now take the function $f$ that is $1$ on the double coset $HqH$ and $0$ on the others.  
The equality
\begin{equation*}
f(q)=\frac1{n}\sum_H f(php\inv)
\end{equation*}
can only be satisfied if all terms $f(php\inv)$ are equal to $1$. This means that $php\inv\in HqH$ for all $h$. In particular for $h=e$. This means that $e\in HqH$ and $HqH=H$. So $q\in H$ and $php\inv\in H$ for all $h$. This holds for all $p$ and therefore $H$ is a normal subgroup.

\ebew

\newpage

\bf The dual $^*$-quantum hypergroup \rm
\nl
There are two possible approaches in line with the two options discussed before. We can use the general construction of the dual of $A$ as in the previous section, or we can construct another finite $^*$-quantum hypergroup $B$, then define the pairing and argue that the new quantum hypergroup can be identified with the dual of the first one. We will use the second method as it is more instructive.
\ssnl
However, still we have to make a choice. We can either take the direct construction as it could be done also for the previous case, or we can construct the dual $B$ from the group algebra $B_0$, as we did for $A$. We will follow this path.
\ssnl
For this purpose, we use the following element in the group algebra $\mathbb C G$ of $G$, as before denoted by $B_0$.

\prop
Define $u\in B_0$ by $u=\frac1{n}\sum_{h\in H}\lambda_h$. Then $u$ is a self-adjoint idempotent satisfying
\begin{equation*}
\Delta_0(u)(1\ot u)=\Delta_0(u)(u\ot 1)=u\ot u.
\end{equation*}
\eprop
\bew
It is easy to see that $u^2=u$ and $u^*=u$. Further
\begin{equation*}
\Delta_0(u)(1\ot u)=\frac1{n^2}\sum_{h,k\in H}\lambda_h\ot \lambda_{hk}=\frac1{n^2}\sum_{h,k\in H}\lambda_h\ot\lambda_k=u\ot u.
\end{equation*}
Similarly $\Delta_0(u)(u\ot 1)=u\ot u$.
\ebew

So $u$ is a \emph{group-like} projection as defined in \cite{La-VD1a}.
\ssnl
We can now construct the finite $^*$-quantum hypergroup $B$.

\stel\label{stel:4.7}
Let $B$ be the subalgebra of $B_0$ of elements $ubu$ with $b\in B_0$. Define $\Delta$ on $B$ by $\Delta(b)=(u\ot u)\Delta_0(b)(u\ot u)$. Then $(B,\Delta)$ is a finite $^*$-quantum hypergroup.
\estel

\bew
i) It is clear that $\Delta$ is a linear map from $B$ to $B\ot B$. To show that it is coassociative we use that $\Delta_0$ is coassociative and that 
\begin{equation*}
\Delta(b)=(u\ot 1)\Delta_0(b)(u\ot 1)=(1\ot u)\Delta_0(b)(1\ot u)
\end{equation*}
when $b\in B$. This follows from the fact that 
\begin{align*}
\Delta_0(u)(1\ot u)&=\Delta_0(u)(u\ot 1)=u\ot u\\
(1\ot u)\Delta_0(u)&=(u\ot 1)\Delta_0(u)=u\ot u.
\end{align*}
Also $\Delta$ is unital on $B$ because $u$ is the unit in the subalgebra $B$.
\snl
ii) For the counit $\varepsilon$ on $B_0$ we have, when $b\in B$, 
\begin{equation*}
(\varepsilon\ot\iota)\Delta(b)=u((\varepsilon\ot \iota)\Delta_0(b))(u)=ubu=b.
\end{equation*}
Similarly on the other side. So the restriction of $\varepsilon$ to $B$ is a counit for $\Delta$ on $B$.
\ssnl
iii) For the antipode on $B$ we use the antipode on $B_0$. Because $S(u)=u$, it leaves $B$ invariant and it is an anti-isomorphism also of $B$.
\ssnl
iv) We consider a left integral $\varphi$ on $B_0$, given by $\varphi(\lambda_e)=1$ and $\varphi(\lambda_p)=0$ when $p\neq e$. We restrict it to $B$. This restriction is still faithful.
\ssnl
v) Furthermore 
\begin{equation*}
(\iota\ot \varphi)\Delta(b)=(\iota\ot\varphi)((u\ot 1)\Delta_0(b)(u\ot 1))=u\varphi(b).
\end{equation*}
So $\varphi$ is left invariant on $B$. It is also right invariant.
\snl
And just as before, it follows from the fact that 
\begin{equation*}
S((\iota\ot\varphi)(\Delta_0(b)(1\ot d)))=(\iota\ot\varphi)((1\ot b)\Delta_0(d))
\end{equation*}
for all $b,d$ in the group algebra, that the same formula holds for $\Delta$ on $B$. We again use that $\Delta(b)=(u\ot 1)\Delta_0(b)(u\ot 1)$ for all $b\in B$ and that $S(u)=u$.
\ebew

\opm
Observe the similarity between the  two cases. The map $F:b\mapsto ubu$ is also a conditional expectation from $B_0$ onto $B$. Here however, it is not unital as $F(1)=u\neq 1$. The unit in $B$ is not the unit in $B_0$ while the unit in $A$ is the same as the unit in $A_0$. There is also a different behavior with respect to the counit. For the counit $\varepsilon$ on $B_0$ we do have $\varepsilon\circ F=\varepsilon$ because $\varepsilon(u)=1$. For the conditional expectation $E$ on $A_0$ we do not have $\varepsilon\circ E=\varepsilon$ as we observed already before, see a remark after Proposition \ref{prop:4.2}.
\ssnl
Still the two results are of the same type because $\Delta(b)=(F\ot \iota)\Delta_0(b)=(\iota\ot F)\Delta_0(b)$ for $b\in B$, just like we have for $E$ on $A_0$.
\eopm

When the group $H$ is a normal subgroup, the element $u$ commutes with all elements of $B_0$ and then $\Delta(b)=\Delta_0(b)(1\ot u)$. It follows again that $\Delta$ is a homomorphism on $B$. Remark however that here $\Delta$ is not equal to the restriction of $\Delta_0$ to $A\ot A$  as was the case for $A_0$.
\ssnl
Also the converse is true. Indeed, if $\Delta$ is a homomorphism, then $B$ is a Hopf algebra. We will show below that $B$ is the dual of $A$. Then $A$ is also a Hopf algebra and the coproduct on $A$ is also a homomorphism. We have seen that then $H$ is a normal subgroup of $G$, see Proposition \ref{prop:4.6a}. 
\nl
\bf The pairing of $A$ with $B$ \rm
\nl
Consider the natural pairing of the function algebra $A_0$ with the group algebra $B_0$ of $G$  given by $\langle f,\lambda_p\rangle=f(p)$ for all $p$. We use the conditional expectations $E:A_0\to A$ and $F:B_0\to B$.

\prop
For all $f\in A_0$ and $b\in B_0$ we have $\langle E(f),b\rangle=\langle f, F(b)\rangle$. 
The restriction of the pairing to $A\times B$ is a pairing of $^*$-algebras. The coproducts induced by the pairing are the given coproducts on the algebras $A$ and $B$.
\eprop

\bew
i) Take $f\in C(G)$ and $p\in G$. Then 
\begin{equation*}
\langle E(f),\lambda_p\rangle=E(f)(p)
=\frac1{n^2}\sum_{h,k\in H}f(hpk)
=\frac1{n^2}\sum_{h,k\in H}\langle f, \lambda_{hpk} \rangle
\end{equation*}
and we see that $\langle E(f),\lambda_p\rangle=\langle f, u\lambda_p u\rangle=\langle f, F(\lambda_p)\rangle$.
\ssnl
ii) It follows that the pairing of $A$ with $B$ is still non-degenerate. Indeed, assume  that  $f\in A$ and $\langle f,b\rangle =0$ for all $b\in B$. Then with $b\in B_0$ we have
\begin{equation*}
\langle f,b \rangle=\langle E(f),b\rangle=\langle f,F(b)\rangle=0.
\end{equation*}
Then $f=0$ as the original pairing is non-degenerate. Similarly on the other side.
\ssnl
iii) Next take $f\in A$ and $b,b'\in B$. Then 
\begin{align*}
\langle \Delta(f),b\ot b'\rangle
&=\langle (E\ot\iota)\Delta_0(f),b\ot b'\rangle\\
&=\langle \Delta_0(f),F(b)\ot b'\rangle\\
&=\langle \Delta_0(f),b\ot b'\rangle\\
&=\langle f,bb'\rangle.
\end{align*}
Similarly we have $\langle f\ot f',\Delta(b)\rangle=\langle ff',b\rangle$.
\ssnl
iv) Finally we need to remark that the coproducts are $^*$-maps on both algebras.
\ebew

The construction we have used here, namely to get  quantum hypergroups from quantum groups is found also in \cite{De-VD2}. What we have done here is just a special case of the general procedure used in \cite{De-VD2}. We have another application in Section \ref{s:two-sub}, see Remark \ref{opm:6.17}.

\section{\hspace{-17pt}. Low-dimensional examples} \label{s:two-dim} 

We will now consider a special case of the construction in the previous section. We take for $G$ the group $S_3$ of permutations of the set $\{1,2,3\}$. We denote by $\sigma_1$ and $\sigma_2$ the permutations $(1,2)(3)$ and $(1)(23)$ respectively. They generate the group. They satisfy 
\begin{equation*}
\sigma_1^2=\sigma_2^2=e
\tussenen
\sigma_1\sigma_2\sigma_1=\sigma_2\sigma_1\sigma_2.
\end{equation*}
We denote $\sigma_1\sigma_2\sigma_1$ by $\sigma_3$. It is the permutation $(13)(2)$ and also satisfies $\sigma_3^2=e$. Further let $p=\sigma_1\sigma_2$, then 
\begin{equation*}
p^2=\sigma_1\sigma_2\sigma_1\sigma_2=\sigma_2\sigma_1\sigma_2\sigma_2=\sigma_2\sigma_1.
\end{equation*}
We see that $p^3=e$. This gives the elements of the group: $\{e,p,p^2,\sigma_1,\sigma_2,\sigma_3\}$.
\nl
We now consider the subgroup $H$ with elements $e$ and $\sigma_1$.

\prop
The left cosets are $H=\{e,\sigma_1\}$, $\{p,\sigma_3\}$ and $\{p^2,\sigma_2\}$. The right cosets are $H=\{e,\sigma_1\}$, $\{p,\sigma_2\}$ and $\{p^2,\sigma_3\}$. The double cosets are $H=\{e,\sigma_1\}$ and $V=\{p,p^2,\sigma_2,\sigma_3\}$.
\eprop

This is obtained easily from the properties of the generators $\sigma_1$ and $\sigma_2$. 
\ssnl
The $^*$-algebra $A$ of functions on $S_3$ that are constant on double cosets is spanned by the two  characteristic functions $u=\Chi_H$ and $v=\Chi_V$. They satisfy $u^2=u$, $v^2=v$ and $uv=vu=0$. For the involution we have $u^*=u$ and $v^*=v$. The sum $u+v$ is the identity of $A$.
\ssnl
For the coproduct $\Delta$ on $A$ we get the following from the general theory.

\prop\label{prop:5.2}
We have $\Delta(1)=1\ot 1$ and 
\begin{align*}
\Delta(u)&=u\ot u+\frac12 v\ot v\\
\Delta(v)&=u\ot v+v\ot u+\frac12 v\ot v.
\end{align*}
\eprop

\bew
This is again easy to verify. We have e.g.\
\begin{equation*}
\Delta(u)(p,p)=\frac12 (u(p^2)+u(p\sigma_1 p))=\frac12
\end{equation*}
because $p^2\notin H$ while $p\sigma_1 p=\sigma_1\in H$. Besides we get
\begin{equation*}
\Delta(u)(e,e)=u(e)=1
\tussenen 
\Delta(u)(e,p)=\Delta(u)(p,e)=0
\end{equation*}
because $p\notin H$. This proves the formula for $\Delta(u)$. The formula for $\Delta(v)$ can be obtained from this and $\Delta(1)=1\ot 1$.
\ebew

For the counit, we get $\varepsilon(u)=1$ and $\varepsilon(v)=0$. For the antipode we have $S(u)=u$ and $S(v)=v$. \ssnl
Let us finally look at the integrals. Because the coproduct is coabelian, the left integral is also a right integral. We obtain it from the general formula.

\prop
The left integral $\varphi$ is given by $\varphi(u)=2$ and $\varphi(v)=4$.
\eprop

We verify the necessary formulas. We have
\begin{align*}
(\iota\ot\varphi)\Delta(u)&=\varphi(u)u+\frac12\varphi(v)v=2(u+v)=\varphi(u)1\\
(\iota\ot\varphi)\Delta(v)&=\varphi(v)u+\varphi(u)v+\frac12\varphi(v)v=4(u+v)=\varphi(v)1.
\end{align*}
This shows that indeed $\varphi$ is left invariant.
\ssnl
In this case, it follows that it is a left integral. Indeed, first we have
\begin{equation}
(\iota\ot\varphi)(\Delta(x)(1\ot y))=(\iota\ot \varphi)((1\ot x)\Delta(y))\label{eqn:5.1}
\end{equation}
for all $x$ and $y$ in $A$. 
When $x=1$ or $y=1$,  (\ref{eqn:5.1}) follows from the left invariance of $\varphi$. On the other hand, when $x=u$ and $y=u$, it follows because the algebra is abelian. Then (\ref{eqn:5.1}) is true for all $x,y$ as the algebra is spanned by $1$ and $u$. This proves that we have a left integral whose antipode $S$ is the identity map.
\ssnl
In fact, as we see, a left invariant functional on such a two-dimensional abelian algebra is always a left integral for the trivial antipode. 
\nl
All of this suggest a set of examples treated in the following propositions.

\prop\label{prop:5.4} 
Let $A$ be the algebra spanned by the identity $1$ and an idempotent element $v$. Take any complex number $\alpha$ not equal to $-1$. Define $\Delta$ on $A$ by $\Delta(1)=1\ot 1$ and
\begin{equation*}
\Delta(v)=v\ot 1+ 1\ot v+ \alpha\, v\ot v.
\end{equation*}
Then $(A,\Delta)$ is a quantum hypergroup. If $\alpha$ is real and $v^*=v$, then it is a $^*$-quantum hypergroup.
\eprop
\bew
i) It is easy to verify that $\Delta$ is coassociative. If we let $\varepsilon(1)=1$ and $\varepsilon(v)=0$, we see that $\varepsilon$ is a counit for $\Delta$. It is a homomorphism. 
\ssnl
ii) Define $\varphi$ by $\varphi(v)=1$ and $\varphi(1)=-\alpha$. Then
\begin{equation*}
(\iota\ot\varphi)\Delta(v)=\varphi(1)v+\varphi(v)1+\alpha\varphi(v)v=\varphi(v)1.
\end{equation*}
Because we always have that $(\iota\ot\varphi)\Delta(1)=\varphi(1)1$ when $\Delta(1)=1\ot 1$, we see that $\varphi$ is left invariant. 
\ssnl
iii) If $u=1-v$ we have $\varphi(u)=\varphi(1)-\varphi(v)=-(1+\alpha)$. Because we assume that $\alpha\neq -1$, we see that $\varphi$ is faithful.
\ssnl
iv) That $\varphi$ is a left integral for the antipode defined as the identity follows as above. Remark that the identity map is an anti-isomorphism because the algebra is abelian.
\ssnl
v) The last statement is obvious. Remark that $a\mapsto S(a)^*$ is the involution itself and it is a homomorphism because the algebra is abelian.
\ebew

For $\Delta(u)$ and $\Delta(v)$  we have
\begin{equation*}
\Delta(u)=u\ot u -(\alpha+1)\, v\ot v
\qquad\text{and}\qquad
\Delta(v)=u\ot v+v\ot u+(\alpha+2)v\ot v.
\end{equation*}
We see that for $\alpha=-1$ we have $\Delta(u)=u\ot u$. 
Then for a left integral we need $\varphi(u)=0$ and $\varphi$ can not be faithful.
\ssnl
For $\alpha=-2$ we get $\Delta(u)=u\ot u+v\ot v$ and $\Delta(v)=u\ot v+v\ot u$. Only in this case $\Delta$ is a homomorphism. This gives the Hopf algebra of the group with two elements.
\ssnl
For $\alpha=-\frac32$ we get the first example as in Proposition \ref{prop:5.2}.
\nl
For the dual of this more general example we have to consider two cases.

\prop\label{prop:5.5}
Let $A$ be the algebra generated by $1$ and an element $v$ satisfying $v^2=v$ and let $B$ be the algebra generated by $1$ and an element $w$ satisfying $w^2=w$. Assume that $\alpha\neq 0$. Define a pairing of $A$ with $B$
 by
\begin{align*}
\langle 1,1\rangle =1 &\tussenen \langle v,1 \rangle =0 \\
\langle 1,w\rangle=0 &\tussenen \langle v,w\rangle =\alpha\inv.
\end{align*}
Then the coproduct $\Delta$ on $A$ is unital and $
\Delta(v)=v\ot 1 + 1\ot v + \alpha\, v\ot v$
\eprop
\bew
Because $(1,v')$ with $v'=\alpha v$, and $(1,w)$ are dual bases of $A$ and $B$ respectively, the coproduct can easily be derived from the product. Indeed, because in $B$ we have 
\begin{align*}
1\cdot 1=1 &\tussenen 1\cdot w=w\\
w\cdot 1=w &\tussenen w\cdot w=w,
\end{align*}
we must have 
\begin{align*}
\Delta(1)&=1\ot 1 \\
\Delta(v')&=1\ot v'+ v'\ot 1+v'\ot v'.
\end{align*}
If we replace $v'$ by $\alpha v$ in the last equality, we get the result.
\ebew

By symmetry, the coproduct on $B$ is also unital and $\Delta(w)=1\ot w + w\ot 1 + \alpha w \ot w$.
\ssnl
We get the following duality result.

\prop
If $\alpha\neq -1$ and $\alpha\neq 0$ then the above pairing is a pairing of quantum hypergroups. 
\eprop

Indeed, the pairing of $A$ with $B$ induces the coproduct as in Proposition \ref{prop:5.4}. Hence $A$ is a quantum hypergroup. We need $\alpha\neq -1$ for this to hold. And we need $\alpha\neq 0$ for Proposition \ref{prop:5.5} to hold.
\ssnl
However, we also have the following result.

\prop
Let $A$ be generated by $1$ and $v$ satisfying $v^2=v$. Now let $D$ be generated by $1$ and an element $y$ satisfying $y^2=0$. Again let $\alpha\neq 0$. Define a pairing by
\begin{align*}
\langle 1,1\rangle =1 &\tussenen \langle v,1 \rangle =0 \\
\langle 1,y\rangle=0 &\tussenen \langle v,y\rangle =\alpha\inv.
\end{align*}
The induced coproduct $\Delta$  on $A$ is unital and satisfies
\begin{equation*}
\Delta(v)=v\ot 1+1\ot v
\end{equation*}
while the coproduct on $D$ is unital and 
\begin{equation*}
\Delta(y)=y\ot 1+1\ot y+\alpha y\ot y.
\end{equation*}
\eprop
\bew 
Again we have the dual bases $(1,\alpha v)$ and $(1,y)$ and we can derive the coproduct on $A$ easily from the product in $B$ as in the proof of the Proposition \ref{prop:5.5}. Because $1$ is the unit and $y^2=0$ we now get
$\Delta(1)=1$ and
\begin{equation*}
\Delta(v)=1\ot v+v\ot 1.
\end{equation*}
The coproduct on $D$ is the same as before because the product on $A$ is the same. So $\Delta$ on $D$ is unital and
\begin{equation*}
\Delta(y)=1\ot y+y\ot 1+\alpha y\ot y.
\end{equation*}

\vskip -20pt
\ebew

Remark that we can assume that $\alpha=1$ by replacing $y$ by $\alpha y$.
\ssnl
We get again a pairing of quantum hypergroups.

\prop
The algebra $D$ with the coproduct as in the above proposition is again a quantum hypergroup. As in the previous proposition, we have a pairing of quantum hypergroups.
\eprop

The left integral on $A$ is given by $\varphi(v)=1$ and $\varphi(1)=0$. The left integral on $D$ is as before, we have $\varphi(y)=1$ and $\varphi(1)=-\alpha$.

\ssnl
This suggest still another possibility.

\prop Let $C$ be the unital algebra generated by an element $x$ satisfying $x^2=0 $ and $D$ the unital algebra generated by $y$ satisfying $y^2=0$. Again let $\alpha\neq 0$. Define a pairing by 
\begin{align*}
\langle 1,1\rangle =1 &\tussenen \langle x,1 \rangle =0 \\
\langle 1,y\rangle=0 &\tussenen \langle y,y\rangle =\alpha\inv.
\end{align*}
Then the coproducts on $C$ and on $D$ are  unital and further given by the same formula:
\begin{align*}
\Delta(x)&=x\ot 1+ 1\ot x \\
\Delta(y)&=y\ot 1+ 1\ot y.
\end{align*}
\eprop

Again we can assume that $\alpha=1$. And we get a pairing of quantum hypergroups.

\opm
In all these cases we have a pairing of $^*$-quantum hypergroups. The generators $v,w,x,y$ are assumed to be self-adjoint. Only for the first case (Proposition \ref{prop:5.5}), we have a good involutive structure in the sense that we have operator algebras. Indeed, if $\varphi(1)=0$, then $\varphi$ can not be positive. In an operator algebra, we also can not have non-zero self-adjoint elements $y$ satisfying $y^2=0$.
\eopm
\snl
\bf Simple examples with no invariant functionals or no integrals \rm
\nl
We will now consider some examples, inspired by the previous ones, to illustrate the difference between the existence of invariant functionals and integrals. 
\ssnl
The first example  is derived from a groupoid with only two elements. This is an example  where there are integrals but no invariant functionals.

\voorb\label{voorb:5.11}
Take the groupoid with two elements $x,y$ where only the products $xx$ and $yy$ are defined and satisfy $xx=x$ and $yy=y$.  For the algebra $A$ we take the function algebra. In this case it is $\mathbb C^2$. If $p,q$ are the $\delta$-functions in the points $x,y$ respectively, the coproduct is given by
\begin{equation*}
\Delta(p)=p\ot p \tussenen \Delta(q)=q\ot q.
\end{equation*}
The coproduct is not unital. As $1=p+q$ we have  $\Delta(1)=p\ot p+q\ot q$ and this is not $1\ot 1$. For the counit $\varepsilon$ we have $\varepsilon(p)=\varepsilon(q)=1$. This is not a homomorphism because $pq=0$ so that $\varphi(pq)\neq \varepsilon(p)\varepsilon(q)$. 
\ssnl
We can not have non-zero invariant functionals. On the other hand, any linear functional is a left integral for the antipode equal to the identity map. Indeed, we claim that 
\begin{equation*}
(\iota\ot\varphi)(\Delta(a)(1\ot b))=(\iota\ot\varphi)((1\ot a)\Delta(b))
\end{equation*}
for all $a,b$. To show this, assume that $a,b$ are $p$ or $q$. Then
\begin{equation*}
(\iota\ot\varphi)(\Delta(a)(1\ot b))=a\varphi(ab) 
\tussenen
(\iota\ot\varphi)((1\ot a)\Delta(b))=b\varphi(ba).
\end{equation*}
If $a\neq b$ we have $ab=ba=0$ and $\varphi(ab)=\varphi(ba)=0$ so that we have equality in all the cases.
\ssnl
The dual algebra is the groupoid algebra spanned by the elements $x,y$. It is again $\mathbb C^2$ and the pairing is given by 
\begin{align*}
\langle p,x\rangle =1 &\tussenen \langle p,y \rangle =0 \\
\langle q,x\rangle=0 &\tussenen \langle q,y\rangle =1.
\end{align*}
\evoorb

Next we look for examples of finite-dimensional algebras where we have invariant functionals but no integrals. We will eventually arrive at a good example, but we find it via different steps.  In all our cases, we have the identity $1$ in $A$ and the integral on $A$ given by a cointegral $h$ in $B$. 

\voorb
First we try the case  $B=\mathbb C^2$. The algebra $B$ is spanned by orthogonal idempotents $h,x$. The algebra $A$ is spanned by elements $1,u$ where $1$ is the identity. For the pairing we take
\begin{align*}
\langle 1,h\rangle =1 &\tussenen \langle 1,x \rangle =0 \\
\langle u,h\rangle=0 &\tussenen \langle u,x\rangle =1.
\end{align*}
Then the coproduct on $A$ is given by
\begin{equation*}
\Delta(1)=1\ot 1 \tussenen \Delta(u)=u\ot u.
\end{equation*}
A left invariant functional $\varphi$ on $A$ is given by $\varphi(a)=\langle a,h\rangle$. 
\ssnl
We claim that $\varphi$ is always a left integral with antipode equal to the identity map. For this we need to have
\begin{equation*}
(\iota\ot\varphi)(\Delta(a)(1\ot b))=(\iota\ot\varphi)((1\ot a)\Delta(b))
\end{equation*}
for all $a,b$. This is true when either $a$ of $b$ is equal to $1$ because $\Delta(1)=1\ot 1$ and $\varphi$ is an invariant functional. On the other hand it is also true for $a=b=u$ because $\Delta(u)=u\ot u$.
\evoorb

In fact, we have seen already that an invariant functional on a two-dimensional abelian algebra is always an integral, see a remark before Proposition \ref{prop:5.4}.
\ssnl
We push the previous example a little further by looking at $\mathbb C^3$.

\voorb\label{voorb:5.13}
We take for both $A$ and $B$ the commutative algebra $\mathbb C^3$. We use $e_1,e_2,e_3$ for the orthogonal idempotents in $A$ and $h,x,y$ for the orthogonal idempotents in $B$. The identity in $A$ is $e_1+e_2+e_3$. We consider two other elements $u$ and $v$ given by
\begin{equation*}
u=e_1+e_2-e_3 \tussenen v=e_1-e_2+e_3.
\end{equation*}
We have $u^2=v^2=1$. 
The elements $(1, u,v)$ is  still a basis in $A$ and we can define a pairing of $A$ with $B$ so that this is a dual basis of $(h,x,y)$. So 
\begin{align*}
\langle 1, h \rangle&=1  & \langle u,h \rangle&=0  &  \langle v,h \rangle&=0 \\
\langle 1, x \rangle&=0  & \langle u,x \rangle&=1 &  \langle v,x \rangle&=0 \\
\langle 1, y \rangle&=0  & \langle u,y \rangle&=0  &  \langle v,y \rangle&=1.
\end{align*} 
The coproduct on $A$ induced by this pairing is given by\begin{equation*}
\Delta(1)=1\ot 1 \qquad
\Delta(u)=u\ot u \qquad
\Delta(v)=v\ot v.
\end{equation*}
We define $\varphi$ on $A$ by $\varphi(a)=\langle a,h\rangle$. So $\varphi(1)=1$ and $\varphi(u)=\varphi(v)=0$.
In this case it follows that $\varphi(e_1)=0$ and $\varphi(e_2)=\varphi(e_3)=\frac12$.
\ssnl
 Then $\varphi$ is left invariant:
\begin{align*}
(\iota\ot\varphi)\Delta(1)&=\varphi(1) 1 \\
(\iota\ot\varphi)\Delta(u)&=\varphi(u) u=0=\varphi(u)1 \\
(\iota\ot\varphi)\Delta(v)&=\varphi(v) v=0 =\varphi(v)1.
\end{align*}
It is also right invariant with the same argument, but also because $B$ is abelian. It is the unique invariant functional because $\varphi(1)\neq 0$.
\ssnl
Suppose that $\varphi$ is a left integral with antipode  $S$.  Then by definition we have 
\begin{equation*}
S((\iota\ot\varphi)(\Delta(a)(1\ot c)))=(\iota\ot\varphi)((1\ot a)\Delta(c).
\end{equation*}
If we take for $a$ or $c$ any of the elements $1$, $u$, $v$ this reads as
\begin{equation*}
S(a)\varphi(ac)=c\varphi(ac).
\end{equation*}
Because $a^2=1$ for these $3$ elements, and because $\varphi(1)=1$, we must have $S=\iota$. But with $a=u$ and $c=v$ we find $\varphi(uv)=\varphi(e_1-e_2-e_3)=-1$. We get a contraction as $\varphi(uv)\neq 0$ and $u\neq v$.
\ssnl
We see that $\varphi$ is invariant but not an integral. Unfortunately, this integral is not faithful because 
$\varphi(e_1)=0$.
\evoorb

In order to get an example with a faithful left invariant functional that is not a left integral, we try the above example in $\mathbb C^4$.

\voorb
i) We take for $A$ and $B$ the algebra $\mathbb C^4$. We take the basis of orthogonal idempotents in $A$ and $B$ and denote them by $e_1,e_2,e_3,e_4$ in $A$ and $h,x,y,z$ in $B$. We take
\begin{align*}
u&=e_1+e_2-e_3-e_4 \\
v&=e_1-e_2-e_3+e_4 \\
w&=e_1-e_2+e_3-e_4 \\
1&=e_1+e_2+e_3+e_4.
\end{align*}
These elements still form a basis and we can define the pairing such that this is a dual basis as in the previous example.
As before we take $\varphi(a)=\langle a,h\rangle$ for all $a$. We get $\varphi(u)=\varphi(v)=\varphi(w)=0$ while $\varphi(1)=1$. This implies $\varphi(e_i)=\frac14$ for all $i$. Again $\varphi$ will be left and right invariant and because $\varphi(1)=1$ it is the unique invariant functional.
\ssnl
Now suppose that it is a left integral with antipode $S$. As in the previous example, we have $S=\iota$ because $u^2=v^2=w^2=1$. 
\ssnl
As we have
\begin{equation*}
uv=w\tussenen vw=u \tussenen wu=v,
\end{equation*}
 now $\varphi(uv)=\varphi(vw)=\varphi(uw)=0$. Therefore, we do get a left integral. And it is faithful.
 \ssnl
ii) Fortunately, it is easy to modify the example. We replace $w$ by $w':=w+\lambda 1$  for some $\lambda\in \mathbb C$. We will still get a basis for $A$ and we can again take $1,u,v,w'$ as a dual basis to define the pairing. The functional $\varphi$ is invariant. If it is a left integral with antipode $S$, we still will have $S(u)=u$ and $S(v)=v$ because $u^2=v^2=1$. Because $uv=w$ we get $\varphi(uv)=-\lambda$ because $\varphi(w')=0$. If $\lambda\neq 0$ we arrive at a contradiction.
\evoorb

We have a similar situation in the following example.

\voorb
i) We again take for $B$ the algebra $\mathbb C^4$, but now we take for $A$ the algebra $M_2(\mathbb C)$ of $2\times 2$ matrices over $\mathbb C$. In $B$ we use the basis $(h,x,y,z)$ of mutually orthogonal idempotents  like in the previous example. In $A$ we use the matrix elements $(e_{ij})$. We define a pairing by defining elements $u,v,w$ in $A$ such that $(1,u,v,w)$ is a dual basis in $A$ for the basis we have in $B$. Here $1$ is the identity in $A$.
\ssnl
ii) For the coproduct $\Delta$ on $A$ we have $\Delta(a)=a\ot a$ for $a$ equal to $1, u, v$ or $w$  as we have in the previous examples. The functional $\varphi$ on $A$ defined by $\varphi(a)=\langle a,h\rangle$ is left and right invariant.
\ssnl
iii) We now take $u=e_{12}+e_{21}$ and $v=e_{12}-e_{21}$. Then $u^2=-v^2=1$. For $w$ we take $pe_{11}+qe_{22}$ where $p,q\in \mathbb R$. Suppose that $\varphi$ is a left integral with antipode $S$. This means that $S(a)\varphi(ac)=c\varphi(ac)$ for $a$ and $c$ any of the elements $1,u,v,w$ as we have seen in  Example \ref{voorb:5.13}. With $a=c=u$ we find $S(u)=u$ because $\varphi(u^2)=\varphi(1)=1$ and with $a=c=v$ we find $S(v)=v$  because $\varphi(v^2)=-\varphi(1)=-1$.
\ssnl
iv) Now we want to have that $\varphi(uv)\neq 0$, then we get a contradiction because it would follow that $S(u)=v$ and because $S(u)=u$ this is a contradiction. Now $uv=e_{22}-e_{11}$. Then $\varphi(uv)=0$ if $\varphi(e_{11)}=\varphi(e_{22})$ We know that $\varphi(w)=0$ and this would imply then that $p+q=1$.
\ssnl
v) It follows that we have to take $p\neq q$ and $p\neq -q$. Finally, we want a faithful functional, so we must also exclude the case where $\varphi(e_{11})=0$. However, if that is the case, from $\varphi(w)=0$ we get $q\varphi(e_{22})=0$. If $q\neq 0$ we get $\varphi(e_{22})=0$ and $\varphi(1)=1$. Therefore we take $q\neq 0$. Similarly we need to take $p\neq 0$. 
\ssnl
vi) So if $w=pe_{11}+qe_{22}$ where 
\begin{equation*}
p\neq q \qquad p\neq -q \qquad p\neq 0 \qquad q\neq 0,
\end{equation*}
we get an example of a pairing where $A$ has a faithful invariant functional that is not a left integral.
\evoorb

\section{\hspace{-17pt}. Quantum hypergroups from  two finite subgroups} \label{s:two-sub}

In this section, we consider some examples of finite quantum hypergroups of a different nature. 
They are inspired by the theory of bicrossproducts. We will formulate the results, but only for a few we will include a proof. Most of the arguments are obvious and easy, so we leave them as an exercise to the reader.  Moreover, details can be found for  the more general cases in \cite{La-VD4}
\nl
\bf Finite quantum hypergroups from a  group with two finite subgroups \rm
\nl
The starting point is any group $G$, not necessarily finite, with two finite subgroups $H$ and $K$. The only requirement is that  $H\cap K=\{e\}$ where $e$ is the identity of $G$. 

\notat
Denote by $\Omega$ the set of pairs $(h,k)$ in $H\times K$ satisfying $hk\in KH$. For $(h,k)\in \Omega$ we define elements $h\tr k\in K$ and $h\tl k\in H$ by $hk=(h\tr k)(h\tl k)$. These elements are well-defined because $hk\in KH$ and $H\cap K=\{e\}$.
\enotat

One easily shows that $(h,k)\mapsto h\tr k$ is a partially defined left action of the group $H$ on the set $K$ and that $(h,k)\mapsto h\tl k$ is a partial right action of the group $K$ on the set $H$. The actions are unital. On the other hand, for all $h,h'\in H$ and $k,k'\in K$ we have $h\tr e=e$, $e\tl k=e$,
\begin{equation*}
h\tr kk'=(h\tr k)((h \tl  k)\tr k')
\tussenen
h'h\tl k=(h'\tl (h\tr k))(h\tl k)
\end{equation*}
whenever these expressions are defined. 
\ssnl
All these formulas follow easily from the definitions.
\nl
We  introduce two finite-dimensional unital $^*$-algebras $A$ and $B$. The underlying space for both  is the set $C(\Omega)$ of complex functions on $\Omega$, but the product and the involutions are different. 

\prop\label{prop:6.2}
Let $f, f_1,f_2\in C(\Omega)$. If we define $f_1f_2$ and $f^*$ in $C(\Omega)$ by
\begin{align}
(f_1f_2)(h,k)&=\sum_v f_1(h,v)f_2(h\tl v, v\inv k)\label{eqn:6.1}\\
f^*(h,k)&=\overline{f(h\tl k,k\inv)}\nonumber,
\end{align}
then $C(\Omega)$ is a unital $^*$-algebra. The sum is taken over the elements $v$ in $K$ satisfying $(h,v)\in\Omega$. The unit is the function $f$ given by $f(h,k)=0$ except for $k=e$ and $f(h,e)=1$ for all $h$.
\eprop

Having a sum like in Equation (\ref{eqn:6.1}), to be precise, we always should mention explicitly the set over which the sum is taken. Most of the time however, this should be clear. If that is not the case, we will indicate it. 

\prop\label{prop:6.3}
Let $g,g_1,g_2\in C(\Omega)$. If we define $g_1g_2$ and $g^*$ in $C(\Omega)$ by
\begin{align*}
(g_1g_2)(h,k)&=\sum_u g_1(hu\inv,u\tr k)g_2(u, k)\\
g^*(h,k)&=\overline{g(h\inv,h\tr k)},
\end{align*}
then $C(\Omega)$ is a unital $^*$-algebra. We take the sum over elements $u\in H$ satisfying $(u,k)\in\Omega$. The unit is the function $g$ given by $g(h,k)=0$ except for $h=e$ and $g(e,k)=1$ for all $k$.

\eprop
These two algebras are groupoid algebras, see Remark \ref{opm:6.14} below.
\ssnl
We will use $C(\widehat\Omega)$ when we consider $C(\Omega)$ with this second product and involution. We keep the notation $C(\Omega)$ for the first algebra. We will systematically use the \emph{algebra} $C(\Omega)$ in the first case and the \emph{algebra} $C(\widehat\Omega)$ in the second case. If we only consider the vector space, we simply mention it as the space $C(\Omega)$. Besides, we will systematically use $f,f_1,f_2$ for elements in the algebra $C(\Omega)$ and $g,g_1,g_2$ for elements in the algebra $C(\widehat\Omega)$.
\snl
We have the following result.

\stel\label{stel:4.24}
Define a pairing of the algebras $C(\Omega)$ and $C(\widehat\Omega)$ by 
$$\langle f,g\rangle=\sum_{(h,k)\in\Omega}f(h,k)g(h,k).$$ 
 Then we have a dual pair of finite $^*$-quantum hypergroups. If $HK=KH$ we have a dual pair of Hopf $^*$-algebras.
\estel

We split the proof in a couple of partial results.
\ssnl
First we have the necessary property of the counits, induced by the pairing.

\prop
The counits are given by
\begin{align*}
\varepsilon(f)&=\sum_{k\in K} f(e,k)\quad \text{ for } f\in  C(\Omega),\\
\varepsilon(g)&=\sum_{h\in H} g(h,e)\quad \text{ for } g\in C(\widehat\Omega).
\end{align*}
The counits are $^*$-homomorphisms.
\eprop

For the coproducts, we get the following formulas.

\prop
The coproducts are given by
\begin{align*}
\Delta(f)(u,v;h,k)&= f(uh,k)\delta_{h\tr k}(v) \quad \text{ for } f\in  C(\Omega),\\
\Delta(g)(u,v;h,k)&=g(u,vk)\delta_{u\tl v}(h) \quad\text{ for } g\in C(\widehat\Omega).
\end{align*}
where $(u,v)$ and $(h,k)$ are $\Omega$. They are $^*$-maps and they are unital. 
\eprop
 We use $\delta$ for the Kronecker delta. 
 \ssnl
In general, the coproducts are not homomorphisms. There are simple examples to illustrate this. In fact, we have the following result, see Proposition 3.12 in \cite{La-VD4}.

\prop
The coproducts are homomorphisms if and only if $HK=KH$.
\eprop

\bew
i) Let $f_1,f_2\in C(\Omega)$. Then
\begin{align*}
\Delta(f_1f_2)(u,v;h,k)
&=(f_1f_2)(uh,k)\delta_{h\tr k}(v) \\
&=\sum_r f_1(uh,r)f_2(uh\tl r,r\inv k)\delta_{h\tr k}(v).
\end{align*}
The sum is taken over elements $r\in K$ satisfying $(uh,r)\in\Omega$.  Then $uh\tl r$ is defined.
We also have $(h,k)\in\Omega$ because  we define $\Delta(f_1f_2)$ on $\Omega\times\Omega$.
\ssnl
On the other hand
\begin{align*}
(\Delta(f_1)\Delta(f_2))&(u,v;h,k)\\
&=\sum_{p,r} \Delta(f_1)(u,p;h,r)\Delta(f_2)(u\tl  p,p\inv v;h\tl r,r\inv k)\\
&=\sum_{p,r} f_1(uh,r)f_2((u\tl p)(h\tl r),r\inv k)\delta_{h\tr r}(p)\delta_{(h\tl r)\tr r\inv k}(p\inv v).
\end{align*}
The  sums are taken over elements $p,r\in K$ satisfying $(u,p),(h,r)\in\Omega$.
\ssnl
When $p=h\tr r$, then $(u\tl p)(h\tl r)=uh\tl r$ and we see that in both formulas, we have the same expressions, except for the conditions imposed by the $\delta$-functions. For the second one we have the sum over $p$ with the conditions that 
\begin{equation*}
p=h\tr r 
\tussenen 
p\inv v=(h\tl r)\tr r\inv k.
\end{equation*}
Then $v=(h\tr r)((h\tl r)\tr r\inv k)=h\tr rr\inv k=h\tr k$.  In the end we find
\begin{align*}
\Delta(f_1f_2)(u,v;h,k)&=\sum_r f_1(uh,r)f_2(uh\tl r,r\inv k)\delta_{h\tr k}(v)\\
(\Delta(f_1)\Delta(f_2))(u,v;h,k)&=\sum_{r,(h,r)\in\Omega}f_1(uh,r)f_2(uh\tl r,r\inv k)\delta_{h\tr k}(v).
\end{align*}
So the real difference between the two expressions is that, in the second case, we have the extra restriction that $(h,r)\in \Omega$. 
If $HK=KH$, so that $\Omega=H\times K$, then the extra condition is automatically fulfilled and we see that $\Delta(f_1f_2)=\Delta(f_1)\Delta(f_2)$.
\ssnl
ii) To prove the converse, assume that we have functions $f_1$ and $f_2$ in $C(\Omega)$. Then we find 
\begin{align*}
\Delta(f_1f_2)(u\inv,e;u,e)
&=\sum_r f_1(e,r)f_2(e\tl r,r\inv e)\delta_{u\tr e}(e)\\
&=\sum_r f_1(e,r)f_2(e,r\inv)
\end{align*}
because $e\tl r=e$ and $u\tr e=e$. For $\Delta(f_1)\Delta(f_2)(u\inv,e;u,e)$ we get the same sum, except that we only sum over elements $r$ with the property that $(u,r)\in \Omega$. Now suppose that $\Omega$ is strictly smaller than $H\times K$. Then we choose 
$(u,v)\notin \Omega$. If now $f_1$ is $1$ in $(e,v)$ and $f_2$ is $1$ in $(e,v\inv)$ and $0$ in other points, then we find
\begin{equation*}
\Delta(f_1f_2)(u\inv,e;u,e)=1
\tussenen
(\Delta(f_1)\Delta(f_2))(u\inv,e;u,e)=0.
\end{equation*}
So $\Delta(f_1f_2)\neq \Delta(f_1)\Delta(f_2)$.
\ebew

\opm
i) If $HK=KH$ then $HK$ is a subgroup and we can as well assume that it is all of $G$. In that case, we have a matched pair of subgroups and we find us in the theory of bicrossproducts (cf.\ \cite{Ma}).
\ssnl
ii) We get an example at the other extreme if $G$ is the free group generated by $H$ and $K$. Then $\Omega$ only contains the pairs $(h,e)$ and $(e,k)$. We will consider this special case further in this section.
\eopm

We can verify that the following linear functionals are invariant.

\prop\label{prop:5.9}
i) Define $\varphi$  on $C(\Omega)$ by
\begin{equation*}
\varphi(f)=\sum_{h\in H} f(h,e).
\end{equation*}
Then $\varphi$ is left and right invariant on $C(\Omega)$ .
\ssnl
ii) Define $\varphi$ on $C(\widehat\Omega)$ by
\begin{equation*}
\varphi(g)=\sum_{k\in K} g(e,k).
\end{equation*}
Then again $\varphi$ is left and right invariant on $C(\widehat\Omega)$.
\eprop

Finally we obtain the antipodes. 

\prop\label{prop:4.31}
i) Define $S$ on $C(\Omega)$ by
\begin{equation*}
S(f)(h,k)=f((h\tl k)\inv,(h\tr k)\inv)
\end{equation*}
for $(h,k)\in \Omega$. Then $S$ is an antipode for the linear functional $\varphi$. 
\ssnl
ii) We have the same formula for the antipode on the dual algebra.
\eprop

Observe that $hk=(h\tr k)(h\tl k)$ so that $(hk)\inv=(h\tl k)\inv(h\tr k)\inv$. If we define $\widetilde f$ on $G$ by $\widetilde f(hk)=f(h,k)$ when $(h,k)\in \Omega$ and $\widetilde f(p)=0$ in other points, we see that 
\begin{equation*}
\widetilde {S(f)}(p)=\widetilde f(p\inv)
\end{equation*}
for all $p\in G$. See also item iv) in Remark \ref{opm:6.14}.
\ssnl
The formulas for the maps $f\to S(f)^*$ on $C(\Omega)$ and $g\mapsto S(g)^*$ on $C(\widehat\Omega)$ are easier to handle. 

\prop
We have
\begin{align*}
S(f)^*(h,k)&=\overline{f(h\inv,h\tr k)} \qquad\text{for } f\in C(\Omega)\\
S(g)^*(h,k)&=\overline{g(h\tl k,k\inv)} \qquad\text{for } g\in C(\widehat\Omega).
\end{align*}
\eprop

\bew
We can verify these formulas easily. Take $f\in C(\Omega)$ and $g\in C(\widehat\Omega)$. Then
\begin{align*}
\langle f,g^*\rangle 
&=\sum f(h,k)g^*(h,k)\\
&=\sum f(h,k)\overline{g(h\inv,h\tr k})\\
&=\sum f(h,h\inv \tr k)\overline{g(h\inv ,k)}\\
&=\sum f(h\inv,h\tr k)\overline{g(h ,k)}.
\end{align*}
All the sums are taken over $(h,k)\in\Omega$. We see that 
$$\langle f,g^*\rangle =\overline{\langle S(f)^*,g\rangle}$$
 if $S(f)^*(h,k)=\overline{f(h\inv,h\tr k)}$. A similar argument works for $S(g)^*$ when $g\in C(\widehat\Omega)$. 
\ebew

We can use this to show that $S$ on $C(\Omega)$ is an antipode for $\varphi$ on $C(\Omega)$.

\prop
Assume that $\varphi$ is the integral on the algebra $C(\Omega)$ as defined in item i) of Proposition \ref{prop:5.9}
Let $f_1,f_2\in C(\Omega)$ and define
\begin{equation*}
f_3=(\iota\ot\varphi)(\Delta(f_1)(1\ot f_2^*))
\tussenen
f_4=(\iota\ot\varphi)(\Delta(f_2)(1\ot f_1^*)).
\end{equation*}
Then $S(f_3)^*=f_4$.
\eprop
\bew
For $f_3$ we find
\begin{align*}
f_3(u,v)
&=\sum_h (\Delta(f_1)(1\ot f_2^*))(u,v;h,e)\\
&=\sum_{h,r} \Delta(f_1)(u,v;h,r)f_2^*(h\tl r,r\inv)\\
&=\sum_h f_1(uh,r)\overline{f_2(h,r)}
\end{align*}
where $r=h\inv \tr v$. The sums are taken over all elements $h\in H$. For the second sum, we assume $(h,r)\in\Omega$. Then 
\begin{align*}
S(f_3)^*(u,v)
&=\overline{f_3(u\inv, u\tr v)}\\
&=\sum_h \overline{f_1(u\inv h,s)}f_2(h,s).
\end{align*}
We sum over elements $h$ satisfying $(h,s)\in \Omega$ where $s=h\inv\tr (u\tr v)=h\inv u\tr v$. If we now replace $h$ by $uh$ we find
\begin{equation*}
S(f_3)^*(u,v)=\sum_h \overline{f_1(h,r)}f_2(uh,r)
\end{equation*}
where again $r=h\inv \tr v$. We see that $S(f_3)^*=f_4$.
\ebew

From this, the basic formula for the antipode follows. The proof is simpler when using the formula for the map $f\mapsto S(f)^*$ than a direct proof using the formula for the antipode as in Proposition \ref{prop:4.31}.
\ssnl
One still has to verify that the antipodes are anti-isomorphisms. Again it is easier to show that the maps $f\mapsto S(f)^*$ and $g\mapsto S(g)^*$ are (conjugate linear) isomorphisms.  The arguments are of the same type as for showing that we have involutive algebras in Propositions \ref{prop:6.2} and \ref{prop:6.3}.
\ssnl
This will  finally complete the proof of the fact that $C(\Omega)$ and $C(\widehat\Omega)$ are finite $^*$-quantum hypergroups.
\ssnl
That they are dual to each other is proven in the following proposition.

\prop
Let $f,f'\in C(\Omega)$ and $g,g'\in C(\widehat\Omega)$. Then
\begin{equation*}
\langle \Delta(f),g\ot g'\rangle =\langle f, gg'\rangle 
\tussenen \langle f\ot f',\Delta(g)\rangle = \langle ff',g\rangle.
\end{equation*}
\eprop

\bew
\begin{align*}
\langle \Delta(f),g\ot g'\rangle
&=\sum_{u,v,h,k}\Delta(f)(u,v;h,k)g(u,v)g'(h,k) \\
&=\sum_{u,v,h,k} f(uh,k)\delta_{h\tr k}(v)g(u,v)g'(h,k) \\
&=\sum_{u,h,k} f(uh,k)g(u,h\tr k)g'(h,k) \\
&=\sum_{u,h,k} f(u,k)g(uh\inv,h\tr k)g'(h,k) \\
&=\sum_{u,k} f(u,k) (gg')(u,k).
\end{align*}
The sums are taken over elements with $(u,v)\in \Omega$ and $(h,k)\in\Omega$. Similarly for the other formula.
\ebew

Also observe that for $f,f'\in C(\Omega)$ 
\begin{equation*}
\varphi(ff')=\sum_h (ff')(h,e)=\sum_{h,v} f(h,v)f'(h\tl v,v\inv).
\end{equation*}
The second sum is over elements $v\in K$ satisfying $(h,v)\in\Omega$. So  the map $f'\mapsto \kappa(f')$ where 
$$\kappa(f')(h,v)=f'(h\tl v,v\inv)$$
 is the bijective map of the space $C(\Omega)$ that realizes the dual of $C(\Omega)$ with $C(\widehat\Omega)$. 
\ssnl
Before we continue with a special case, we want to make the following remark.

\opm\label{opm:6.14}
i) The set $\Omega$ has two groupoid structures. One is given by 
$$(h,k)(h',k')=(h,kk') \qquad\text{if } h\tl k=h'.$$
 The other one is 
$$(h,k)(h',k')=(hh',k') \qquad\text{if } k=h'\tr k'.$$
 The inverse of $(h,k)$ is $(h\tl k,k\inv)$ in the first case and $(h\inv, h\tr k)$ in the second case.
\ssnl
ii) The algebra $A$ is the groupoid algebra of $\Omega$ with the first groupoid structure and $B$ is the groupoid algebra of $\Omega$ with the second one. 
\ssnl
iii) The map $(h,k)\mapsto (h\tl k,k\inv)$ is an anti-isomorphism of the first groupoid and an isomorphism of the second one. This explains why $f\mapsto f^*$ is an anti-isomorphism and $f\mapsto S(f)^*$ a isomorphism of $A$.
\ssnl
iv) Finally remark that 
\begin{equation*}
(hk)\inv=((h\tr k)(h\tl k))\inv=(h\tl k)\inv(h\tr k)\inv.
\end{equation*}
If we put $h'=(h\tl k)\inv$ and $k'=(h\tr k)\inv$ we see that $h'\tr k'=k\inv$ and $h'\tl k'=h\inv$. In other words, if $(h',k')=((h\tl k)\inv,(h\tr k)\inv)$ we have $h'k'=(hk)\inv$. 
\eopm
 
\nl
\bf A special case \rm
\nl
We consider now the case of two finite groups $H$ and $K$ with the property that $HK\cap KH=H\cup K$. 
Then $hk\in KH$ can only occur when either $h=e$ or $k=e$. Hence the set $\Omega$ consist of elements $(h,e)$ and $(e,k)$ where $h\in H$ and $k\in K$. Remark that $(e,e)$ is common to both. This makes the system not completely trivial as we will see.
\ssnl
We have the following characterizations of the algebras $A$ and $B$ in this case.

\prop
Let $A_0$ be the tensor product $F(H)\ot \mathbb C K$ of the function algebra $F(H)$ of the set $H$ with the group algebra $\mathbb C K$ of the group $K$. Then $A$ is the subalgebra of functions $f$ on $H\times K$ with support in $\Omega$. The map $E:A_0\to A$ given by restricting a function to $\Omega$ is a self-adjoint conditional expectation satisfying $E(1)=1$.
\eprop

\bew
i) We identify elements of $A_0$ with complex functions on $H\times K$. The product $ff'$ of functions $f,f'\in A$ is given by the formula
\begin{equation*}
(ff')(h,k)=\sum_v f(h,v)f'(h\tl v,v\inv k).
\end{equation*}
Because $f$ has support in $\Omega$, then either $h=e$ or $v=e$. When $h=e$ we have $h\tl v=e$ and so $h\tl v=h$. This also holds for  $v=e$. It follows that 
\begin{equation*}
(ff')(h,k)=\sum_v f(h,v)f'(h,v\inv k)
\end{equation*}
and we see that $A$ is a subalgebra of $A_0$. 
\ssnl
ii) Now assume that $f$ has support in $\Omega$ and that $f'$ is any function in $A_0$. Consider the product $ff'$ in $A_0$. If we assume that $h\neq e$ then $(ff')(h,k)=f(h,e)f'(h,k)$ and then $ff'(h,e)\neq 0$ implies $f'(h,e)\neq 0$. This means that $E(ff')=fE(f')$. On the other hand, assume that $f'$ has support in $\Omega$. Again when $h\neq e$ we have $(ff')(h,k)=f(h,k)f'(h,e)$ and $E(ff')=E(f)f'$ for the same reason.
\ssnl
iii) The identity in $A_0$ is given by the function $f$ that takes the value $1$ in $(h,e)$ and $0$ in other points of $H\times K$. This has support in $\Omega$ and therefore $E(1)=1$.
\ssnl
iv) Finally, if $(h,k)\in \Omega$, then $(h\tl k,k\inv)\in\Omega$. But for $(h,k)\in \Omega$ we have either $h=e$ or $k=e$ so that $(h\tl k,k\inv)=(h,k\inv)$ in both cases. As the adjoint in $A_0$ is given by $f^*(h,k)=\overline{f(h,k\inv)}$, we see that $A$ is a $^*$-subalgebra of $A_0$ and that $E(f^*)=E(f)^*$. 
\ebew

We have a similar result for the algebra $B_0$, defined as $\mathbb C H\ot F(K)$.
\nl
We now consider the natural pairing of the function algebras with the group algebras and the associated tensor product pairing of $A_0$ with $B_0$, This is a pairing of Hopf $^*$-algebras. 
\ssnl
As the pairing is given by $\langle f,g\rangle=\sum_{h,k}f(h,k)g(h,k)$ we have
$\langle E(f),g\rangle =\langle f,E(g) \rangle$ for all $f\in A_0$ and $g\in B_0$.
\ssnl
It is an immediate consequence that $(E\ot E)\Delta_0(f)=(E\ot \iota)\Delta_0(f)=(\iota\ot E)\Delta_0(f)$ when $f\in A$. Indeed, given any $g,g'\in B_0$ we find
\begin{equation*}
\langle (E\ot\iota)\Delta(f),g\ot g'\rangle=\langle f,E(g)g'\rangle=\langle f,E(E(g)g'\rangle=\langle f,E(g)E(g')\rangle
\end{equation*}
and we see that $(E\ot \iota)\Delta_0(f)=(E\ot E)\Delta_0(f)$. Similarly on the other side.
\ssnl
One can now easily get the main result for this example:

\stel\label{stel:6.16}
The pairing of $A$ with $B$ makes $A$ and $B$ into a dual pair of $^*$-quantum hypergroups.
\estel

\bew
i) The counit $\varepsilon$ on $A$ is the restriction to $A$ of the counit on $A_0$. This follows from the fact that the identity in $B$ is the identity in $B_0$. Similarly for the counit on $B$. Therefore, the counits are homomorphisms. Recall that the counit on $A_0$ is given by $\varepsilon(f)=\sum_k f(e,k)$. 
\ssnl
ii) When $\varphi$ is a left integral on $A_0$ we have, for $f\in A$,
\begin{equation*}
(\iota\ot\varphi)\Delta(f)=(\iota\ot\varphi)(E\ot \iota)\Delta_0(f)=\varphi(f)E(1)=\varphi(f)1.
\end{equation*}
and we see that its restriction to $A$ is left invariant on $A$. Recall that the left integral $\varphi$ on $A_0$ is given by
$\varphi(f)=\sum_h f(h,e)$.
\ssnl
iii) The antipode $S$ of $A_0$ is defined as $S(f)(h,k)=f(h\inv,k\inv)$ and the map $(h,k)\mapsto (h\inv,k\inv)$ clearly leaves $\Omega$ invariant. Then $S$ commutes with $E$. Now we have for all $f,f'\in A$ that 
\begin{equation*}
S((\iota\ot\varphi)(\Delta_0(f)(1\ot f')))=(\iota\ot\varphi)((1\ot f)\Delta_0(f'))
\end{equation*}
and if we apply $E$ we arrive at
\begin{equation*}
S((\iota\ot\varphi)(\Delta(f)(1\ot f')))=(\iota\ot\varphi)((1\ot f)\Delta(f')).
\end{equation*}
Hence the restriction of $\varphi$ to $A$ is a left integral.
\ssnl
It follows that the pair $(A,\Delta)$ is a $^*$-quantum hypergroup and similarly for $(B,\Delta)$. As the coproduct on one algebra comes from the product on the other via the pairing, we find a pair of $^*$-quantum hypergroups. This completes the proof.
\ebew

\opm\label{opm:6.17}
i) When we compare this with the Theorems \ref{stel:4.4} and \ref{stel:4.7} we see that in all these cases we have a conditional expectation $E$ from a bigger algebra $A_0$ to the subalgebra $A$. The algebra $A_0$ is a Hopf algebra. The coproduct $\Delta$ on $A$ is obtained from the coproduct $\Delta_0$ on $A_0$ by the formula $\Delta(a)=(E\ot E)\Delta_0(a)$ when $a\in A$. We use that $E(1)$ is the identity in $A$ when $1$ is the identity in $A_0$. We also use that the antipode $S$ of $A_0$ leaves $A$ globally invariant. Then the crucial property is that 
\begin{equation*}
(E\ot E)\Delta_0(a)=(E\ot\iota)\Delta_0(a)=(\iota\ot E)\Delta_0(a)
\end{equation*}
for $a\in A$. The integrals on $A$ are the restrictions to $A$ of the integrals on $A_0$.
\ssnl
ii) For the dual $B$ we also have such a conditional expectation $F$ on a larger algebra $B_0$. When $A_0$ and $B_0$ are a dual pair of Hopf algebra, we need that 
$$\langle E(a),b\rangle =\langle a,F(b)\rangle$$ 
for $a\in A_0$ and $b\in B_0$. Then the restriction of the pairing to $A\times B$ results in a pairing of quantum hypergroups.
\eopm

The algebra $A_0$ of the previous example is generated by elements $\delta_h\ot \lambda_k$ where $\delta_h$ is the function in $F(H)$ with value $1$ in $h$ and $0$ in other points. Similarly $k\mapsto \lambda_k$ is the canonical embedding of $K$ into the group algebra $\mathbb C K$. The subalgebra $A$ is the span of the elements $\delta_e\ot \lambda_k$ and $\delta_h\ot \lambda_e$. 
\ssnl
We now look at the simplest case where  $H=\{e,h\}$ and $K=\{e,k\}$.
Then the algebra $A$ is abelian and is  spanned by the elements
\begin{equation*}
p=\delta_h\ot \lambda_e,\qquad u=\delta_e\ot \lambda_e,\qquad v=\delta_e\ot\lambda_k.
\end{equation*}\ssnl
We have the following relations
\begin{equation*}
\begin{matrix}
&p^2=p, &\qquad pu=pv=0,\\
&u^2=u, &uv=v,\\
&v^2=u.
\end{matrix}
\end{equation*}
The dual algebra $B$ is spanned by the elements 
\begin{equation*}
v'=\lambda_k\ot \delta_e,\qquad u'=\lambda_e\ot \delta_e,\qquad p'=\lambda_e\ot\delta_k
\end{equation*}
and we have the same relations as above. 
\ssnl
The pairing is given by
\begin{equation*}
\begin{matrix}
&\langle p,v'\rangle=1, &\quad \langle u,v' \rangle =0, &\quad \langle v,v' \rangle=0,\\
&\langle p,u'\rangle=0, &\quad \langle u,u' \rangle =1, &\quad \langle v,u' \rangle=0,\\
&\langle p,p'\rangle=0, &\quad \langle u,p' \rangle =0, &\quad \langle v,p' \rangle=1.
\end{matrix}
\end{equation*}
Since we have dual bases, we can easily find the expressions for the coproduct on $A$ from the product of the generators in $B$. One can also use the general formulas for the coproduct on the function algebras and the group algebras.
\ssnl
Observe that the algebras $A$ and $B$ are both $\mathbb C^3$. In $A$ we can look at 
$$\left\{p,\frac12(u+v),\frac12(u-v)\right\}$$
 and we get $3$ orthogonal projections. The pairing however is not the obvious pairing of $\mathbb C^3$ with itself as can be seen from the pairing of the generators above. If we denote these elements in $A$ by $e_1,e_2,e_3$ and the if we use $f_1,f_2,f_3$ for the elements 
 $$\left\{p',\frac12(u'+v'),\frac12(u'-v')\right\}$$ we find for the pairing 
\begin{equation*}
\begin{matrix}
&\langle e_1,f_1\rangle=0, &\quad \langle e_2,f_1 \rangle =\frac12, &\quad \ \, \langle e_3,f_1\rangle=-\frac12,\\
&\langle e_1,f_2\rangle=\frac12, &\quad \langle e_2,f_2 \rangle =\frac14, &\quad \langle e_3,f_2 \rangle=\frac14,\\
&\quad \langle e_1,f_3\rangle=-\frac12, &\quad \langle e_2,f_3 \rangle =\frac14, &\quad \langle e_3,f_3 \rangle=\frac14.
\end{matrix}
\end{equation*}

\section{\hspace{-17pt}. Conclusions and further research}\label{s:conclusions}  

In \cite{La-VD4} we obtain a pair of $^*$-algebraic quantum hypergroups from a pair $(H,K)$ of subgroups of a group $G$ satisfying $H\cap K=\{e\}$. In \cite{La-VD5} we make an attempt to generalize this result for a pair of closed subgroups of a locally compact quantum group. We no longer get algebraic quantum hypergroups but topological quantum hypergroups. Even in the case where the locally compact group $G$ has a compact open subgroup (see \cite{La-VD7}), contrary to what was expected, we also still do not get a pair algebraic quantum hypergroups.
\ssnl
This has inspired us to initiate the study of what would be considered as topological quantum hypergroups. This is done in \cite{La-VD3b}. In that paper, we emphasize on the development of the notions and suggest further study. We do not need this for the examples we find in  \cite{La-VD4, La-VD5}. 
\ssnl
While writing the paper on topological quantum hypergroups we found it instructive to collect the purely algebraic features  behind this construction by treating finite quantum hypergroups in this paper. Finite quantum hypergroups are  special cases of algebraic quantum hypergroups as studied in \cite{De-VD1,De-VD2}. In this paper, we have a treatment independent of the more general case of algebraic quantum hypergroups, especially written for the purpose of generalizing it to topological quantum hypergroups in \cite{La-VD3b}.
\ssnl
In Section \ref{s:integrals} we have defined left and right integrals and obtained a duality for pairs of finite-dimensional algebras with an antipode and  integrals. Finite quantum hypergroups require one more condition, namely that the coproducts are unital. Duality for finite quantum hypergroups is obtained in Section \ref{s:finite}. In Section \ref{s:two-dim} we constructed examples of pairs where we have invariant functionals but no integrals. It suggests that the existence of integrals is much stronger than the existence of invariant functionals.
\ssnl
The duality of the more general objects,  algebras with an antipode and integrals that are not quantum hypergroups, deserves further study. In particular, there is a need for a suitable name for these objects. In \cite{La-VD3a} we treat such pairs of possibly infinite-dimensional non-degenerate algebras.
 \ssnl
 Finally, we spent a great deal of this paper to illustrate the theory with examples. It would be interesting to look for more examples of finite quantum hypergroups and by doing so understand the algebraic requirements.  



\end{document}